\documentclass[a4legal,11pt]{article}
\usepackage{t1enc}
\usepackage[french]{babel}
\usepackage{amssymb,latexsym,amsfonts}
\usepackage{amsmath,amssymb}
\usepackage{fancyhdr}
\usepackage{url}
\usepackage{longtable,hhline}
\usepackage{graphicx}
\usepackage{graphics,epsfig,color}
\graphicspath{ {images/} }
\usepackage{wrapfig}
\usepackage{eurosym}

\usepackage[latin1]{inputenc}
\usepackage[T1]{fontenc}

\usepackage{amsmath,amssymb,amsbsy}

\newtheorem{thm}{\textsc{Theorem}}
\newtheorem{pro}{\textsc{Proposition}}
\newtheorem{cor}{\textsc{Corollary}}
\newtheorem{lem}{\textsc{Lemma}}
\newtheorem{rmq}{\textsc{Remark}}
\newtheorem{defn}{\textsc{Definition}}
\newtheorem{exemple}{\textsc{Example}}
\newtheorem{app}{\textsc{Application}}

\title{Pluripotential theory on the support of closed positive currents and applications to dynamics in $\mathbb{C}^n$}
\author{Fr\'ed\'eric PROTIN \\ 
\\
INSA de Toulouse, D\'epartement G\'enie Math\'ematique et Mod\'elisation  \\
135 avenue de Rangueil, 31400 Toulouse      \\
{\it{fredprotin@yahoo.fr}}}

\begin{document}
\maketitle

\begin{center}\textbf{Abstract}\end{center}$ $ \\
\textit{We extend certain classical theorems in pluripotential theory to a class of functions defined on the support of a $(1,1)$-closed positive current $T$, analogous to plurisubharmonic functions, called $T$-plurisubharmonic functions.  These functions are defined as limits, on the support of $T$, of sequences of plurisubharmonic functions decreasing on this support. In particular, we show that the poles of such functions are pluripolar sets. We also show that the maximum principle and the Hartogs's theorem remain valid in a weak sense. We study these functions by means of a class of measures, so-called ``pluri-Jensen measures'', about which we prove that they are numerous on the support of $(1,1)$-closed positive currents. We also obtain, for any fat compact set, an expression of its relative Green's function in terms of an infimum of an integral over a set of pluri-Jensen measures. We then deduce, by means of these measures, a characterization of the polynomially convex fat compact sets, as well as a characterization of pluripolar sets, and the fact that the support of a closed positive $(1,1)$-current is nowhere pluri-thin. In the second part of this article, these tools are used to study dynamics of a certain class of automorphisms of $\mathbb{C}^n$ which naturally generalize H\'enon's automorphisms of $\mathbb{C}^2$. First we study the geometry of the support of canonical invariant currents. Then we obtain an equidistribution result for the convergence of pull-back of certain measures towards an ergodic invariant measure, with compact support.}


\section{Introduction}

The main goal of this article is to develop a pluripotential theory for $T$-plurisubhar-monic functions (abbreviated $T$-psh), i.e. functions that are limits, on the support of a $(1,1)$-closed positive current $T$, of sequences of plurisubharmonic functions decreasing on this support. The $T$-psh functions were introduced in \cite{BS} (see also \cite{FaG}), and a more general class of functions, called plurisubharmonic in the direction of the current $T$, is considered in \cite{BS} and \cite{S2} from a different point of view. $T$-psh functions naturally arise in complex dynamics, as in \cite{P}, or in \cite{DS4} and \cite{DDS} where they appear with a different name. The second goal of this article is to apply these results to the study of the dynamics of a certain class of automorphisms of $\mathbb{C}^n$ which generalizes H\'enon's automorphisms of $\mathbb{C}^2$, that we call H\'enon-like automorphisms. The article is structured as follows.

In Section \ref{tpsh}, we develop a pluripotential theory for continuous $T$-psh functions. We proceed as follows.
In Subsection \ref{tpsh1}, we give a more intrinsic characterization of such functions. More precisely, we prove that a continuous function $u:supp(T)\rightarrow \mathbb{R}$ is $T$-psh if and only if it satisfies $\displaystyle u(x)\leq \int_{supp(T)}u \ d\mu_x$ for every $x$ in the support $supp(T)$ of $T$ and every pluri-Jensen measure $\mu_x$ relative to $x$ and with support in $supp(T)$. We have equality in the previous relation when we take the infimum over the set of all such pluri-Jensen measures. We also prove that there are sufficiently many such measures in order to construct a potential theory, analogous to the one for plurisubharmonic functions. We deduce some results, as for example a maximum principle for $T$-psh functions, and the fact that a support of a closed positive $(1,1)$-current is nowhere pluri-thin.

In Subsection \ref{subsectiontwo}, for an open connected bounded hyperconvex set $U\subset\mathbb{C}^n$, we express the value at any $x\in U$ of the Green's function relative to $U$ of a fat compact set $A\subset U$ as the infimum of $-\mu_x(\mathring{A})$ over all pluri-Jensen measures $\mu_x$ relative to $x$ with support in $U$. We deduce a characterization of polynomially convex fat compact sets and of pluripolar sets in terms of pluri-Jensen measures. We further prove that the poles of $T$-psh functions form a pluripolar set.

In Subsection \ref{similar} we pursue the analogy between $T$-psh functions and plurisubharmonic functions. In particular, we obtain a Hartogs-like theorem, an extension property across certain small sets, and a Chern-Levine-Nirenberg inequality. We also prove that, under minor conditions, the limit superior of a locally uniformly bounded from above sequence of $T$-psh functions is a $T$-psh function, except on a pluripolar set.

In Section \ref{dynamics} we apply the previous results to the study of the dynamics of H\'enon-like automorphisms.
In Subsection \ref{3.1} we give an equidistribution lemma for the pull-backs of certain closed currents by a H\'enon-like automorphism, that we further exploit to obtain informations about its Julia set.

In Subsection \ref{3.2}, we give an equidistribution result for the pull-backs of certain measures by a H\'enon-like automorphism. More precisely, we prove that the pull-backs of a probability measure defined on the support of certain currents and which does not charge pluripolar sets, converge in a Ces\`aro sense towards an unique canonical invariant measure. This measure is ergodic, with compact support.

\noindent {\it{Notations.}} We will denote $d := \partial + \overline{\partial}$ the usual exterior derivative and $\displaystyle d^c :=\frac{1}{2 i\pi}(\partial -\overline{\partial})$. Thus $\displaystyle dd^c = \frac{i}{\pi}\partial \overline{\partial}$. These operators will be understood in the sense of currents.\\
\noindent A continuous function $u$ from an open set of $\mathbb{C}^n$ into $\mathbb{R}$ is pluriharmonic (harmonic if $n=1$) if and only if $dd^c u =0$. We denote by $\|S\|:=S\wedge \omega^{n-1}$ the trace measure of a closed positive $(1,1)$-current $S$ on $\mathbb{C}^n$, where $\omega := \frac{1}{2}dd^c \log (1+\|z\|^2)$ is the Fubini-Study form. Then, the mass of the current $S$ is $\|S\|(\mathbb{C}^n)$.\\

\section{ Study of $T$-plurisubharmonic functions}\label{tpsh}

\subsection{Fundamental results on $T$-psh functions}\label{tpsh1}

In the sequel, we denote by $\hat{A}$ the polynomially convex envelope of a compact $A$. Following \cite{FaG}, we define the $T$-psh functions as follows :

\begin{defn}\label{deftpsh} Given a closed positive $(1,1)$-current $T$ defined in an open set $V\subset \mathbb{C}^n$ with support $\mathcal{S}\subset V$, a function $u:\mathcal{S}\rightarrow \mathbb{R}\cup \{-\infty\}$ is called \emph{$T$-plurisubharmonic} (abbreviated \emph{$T$-psh}) if it is locally $\|T\|$-integrable and if there exists an open set $U\subset V$ containing $\mathcal{S}$, and a sequence of plurisubharmonic functions $u_m:U\rightarrow \mathbb{R}$ decreasing on $\mathcal{S}$, such that $(u_m|_{\mathcal{S}})_m$ converges simply towards $u$ on $\mathcal{S}\cap U$. Up to shrinking $U$, we may suppose that the functions $u_m$ are continuous.

More generally, if we replace $\mathcal{S}$ by a compact $K\subset\mathbb{C}^n$ in the previous definition, the function $u:K\rightarrow \mathbb{R}\cup \{-\infty\}$ will be said \emph{$K$-plurisubharmonic}.

The elements of a family $(v_m)_{m\in\mathbb{N}}$ of $T$-psh functions will be said \emph{co-$T$-psh} if for every $m \in\mathbb{N}$, there exists an open set $U\subset V$ containing $\mathcal{S}$, and a sequence of plurisubharmonic functions $(v_{m,i})_i:U\rightarrow \mathbb{R}$ decreasing on $\mathcal{S}\cap U$ towards $v_m$.

\end{defn}

The current $T\wedge dd^c u := dd^c (uT)$ is then well defined, positive and closed (see Lemma 4-13 in \cite{FaG}). \\

The main purpose of this section is to present a more intrinsic characterization of continuous $T$-psh functions. First let us recall some definitions and notations. 

A {\it{pluri-Jensen measure $\mu_x$ relative to $x\in\mathbb{C}^n$}} is a regular probability measure with compact support such that, for any plurisubharmonic function $v$ defined on an open set $U$ containing $x$, we have \begin{equation}\label{defjensen}v(x)\leq \int_U v d\mu_x.\end{equation}
If $v$ is only assumed subharmonic, $\mu_x$ is a {\it{Jensen measure relative to $x$}} (therefore it is also a pluri-Jensen measure). The reader should pay attention to the fact that a pluri-Jensen measure in our terminology is called Jensen measure in \cite{St}, \cite{St2}, and \cite{GR}.

\begin{rmq}\label{remarque1}It is sufficient to suppose that (\ref{defjensen}) holds when $v$ is a continuous plurisubharmonic function in order to ensure that $\mu_x$ is a pluri-Jensen measure, thanks to the monotone convergence theorem, since a plurisubharmonic function is a limit of a decreasing sequence of continuous plurisubharmonic functions in a neighborhood of the support of $\mu_x$. It is even sufficient to suppose that $v$ is of the form $\log{|f|}$ for some holomorphic function $f$, since a classical theorem of Bremermann ensures that a continuous plurisubharmonic function can be approximated locally uniformly by functions of the form $\displaystyle \max_{i=1,...,m}c_i\log{|f_i|}$, $c_i> 0$, where $f_i$ are holomorphic functions.\end{rmq}

Let $B\subset \mathbb{C}^n$. We denote by $PJ_B(x)$ the set of pluri-Jensen measures relative to $x$ with support in $B$. Note that $\delta_x\in PJ_B(x)$ if $x\in B$. Moreover, for a given compact set $A\subset\mathbb{C}^n$, for any $x\in \hat{A}$, there exists a measure $\mu_x\in PJ_{A}(x)$ (this follows from Theorem 1.2.9, Theorem 1.2.14 and Remark 1.2.18 in \cite{St}, or even directly from Theorem 2.2.5 in \cite{St}, see also \cite{DS2}). Note also that an automorphism $Q:\mathbb{C}^n\rightarrow\mathbb{C}^n$ induces a bijection ${Q}_*:PJ_{\mathbb{C}^n}(x)\rightarrow PJ_{\mathbb{C}^n}(Q(x))$.

The following lemma ensures the existence of sufficiently many pluri-Jensen measures on the support of a closed positive current, allowing us to establish a potential theory on it.

\begin{lem}\label{local} Let $T$ be a closed positive $(1,1)$-current in $\mathbb{C}^n$, of support $\mathcal{S}$. Then $\forall r>0$, $\forall x\in \mathcal{S}$, there exists a pluri-Jensen measure relative to $x$ with support in $\partial B(x,r)\cap \mathcal{S}$, where $B(x,r)$ denotes the open ball of center $x$ and radius $r$. \end{lem} 

\noindent \textit{Proof} Let $y\in \partial B(x,r)\cap \mathcal{S}$. Then $y$ is a {\it{peak point}} for the uniform algebra  $\mathcal{A}$ of the functions on $\overline{B}(x,r)\cap \mathcal{S}$ with complex values that are restrictions of holomorphic functions defined in a neighborhood of $\overline{B}(x,r)$ : in other words, there exists $f\in \mathcal{A}$ such that $f(y)=1$ and $|f(z)|<1$ for $\displaystyle z\in \left(B(x,r)\cap \mathcal{S}\right)\setminus \{y\}$.

Indeed, assuming that $x$ is the origin, the restriction of the function $f(z)=\frac{1}{2}(1+<z,\frac{y}{\|y\|}>)$ to $\overline{B}(x,r)\cap \mathcal{S}$ fulfills these conditions. ($<\cdot ,\cdot>$ denotes here the canonical hermitian scalar product on $\mathbb{C}^n$.)

In particular, $y$ is a {\it{strong boundary point}} for $\mathcal{A}$, i.e. for any neighborhood $V$ of $y$ for the topology induced on $\partial B(x,r)\cap \mathcal{S}$, there exists $f\in \mathcal{A}$ such that $\displaystyle \sup_V |f|=f(y)=1$ and $|f(z)|<1$ when $\displaystyle z\in \left(B(x,r)\cap \mathcal{S}\right)\setminus V$.

On the other hand, there is no strong boundary point for $\mathcal{A}$ in $B(x,r)\cap \mathcal{S}$. Indeed, suppose that there exists such a point $z\in B(x,r)\cap \mathcal{S}$. The point $z$ cannot be a peak point, since it follows from Proposition 2.5 in \cite{DL} that the algebra $\mathcal{A}$ cannot have local peak points in $B(x,r)\cap \mathcal{S}$. Thus, for any neighborhood $V$ of $y$, a function $f\in \mathcal{A}$ such that $\displaystyle \sup_V |f|=f(y)=1$ and $|f(z)|<1$ for $\displaystyle z\in \left(B(x,r)\cap \mathcal{S}\right)\setminus \{y\}$ takes the value $1$ on $V\setminus \{z\}$. By the Theorem of isolated zeros applied to the holomorphic function whose restriction is $f$, it follows that $f$ is constant. Consequently, $z$ cannot be a strong boundary point.

Therefore, the strong boundary points of $\mathcal{A}$ belong to $\partial B(x,r)\cap \mathcal{S}$. But by Corollary 7.24 in \cite{St2}, the Shilov boundary $\Gamma$ for $\mathcal{A}$ is the closure of the set of strong boundary points of $\mathcal{A}$, and thus $\Gamma \subset \partial B(x,r)\cap \mathcal{S}$. Then by Lemma 6.2 in \cite{A} and Remark \ref{remarque1}, there exists a pluri-Jensen measure relative to $x$ with support in $\Gamma$ (see also Remark 1.2.18 in \cite{St}).$\blacksquare$\\ 

\begin{rmq}\label{rmq.jensen} Let $T$ be a closed positive $(1,1)$-current in an open subset of $\mathbb{C}^n$, with support $\mathcal{S}$, and $u:\mathcal{S}\rightarrow \mathbb{R}$ a $T$-psh function. The monotone convergence theorem and Lemma \ref{local} imply that $\forall x\in \mathcal{S}$, $\forall \nu_x\in PJ_{\mathcal{S}}(x)$, we have $\displaystyle u(x)\leq \int_{\mathcal{S}}u d\nu_x$. Since $\delta_x\in PJ_{\mathcal{S}}(x)$, we obtain \begin{equation}\label{eq.jensen}\displaystyle u(x)= \inf_{\nu_x\in PJ_{\mathcal{S}}(x)}\int_{\mathcal{S}}u d\nu_x.\end{equation}\end{rmq}

\noindent We deduce a maximum principle for $T$-psh functions which improves Lemma 4.3 in \cite{P} (see also Theorem 3.2 in \cite{S2}):

\begin{cor}\label{maximum} Let $T$ be a closed positive $(1,1)$-current in $\mathbb{C}^n$, of support $\mathcal{S}$, and $u:\mathcal{S}\rightarrow \mathbb{R}$ a $T$-psh function. Let $E\subset\mathcal{S}$ a bounded set, $x\in E$, and $r\in ]0,+\infty[$ such that $E\subset\subset B(x,r)$. Then $$\displaystyle \sup_E u\leq \sup_{\partial B(x,r)\cap \mathcal{S}}u.$$In particular, $u$ does not have a strict local maximum.\end{cor}

\noindent \textit{Proof} Let $x\in E$. By Lemma \ref{local}, we can find a pluri-Jensen measure $\mu_x$ with support in $\partial B(x,r)\cap\mathcal{S}$. Then Remark \ref{rmq.jensen} gives $\displaystyle u(x)\leq \int_{B(x,r)\cap\mathcal{S}}u d\mu_x\leq \sup_{\partial B(x,r)\cap\mathcal{S}}u$.$\blacksquare$\\ 

Note that, by virtue of Proposition 4.8.2 in \cite{K}, the previous lemma implies that the support of a closed positive $(1,1)$-current is nowhere {\it{pluri-thin}}, in the terminology of \cite{K}. 

For a subset $A \in \mathbb{C}^n$ and a function $w: A\rightarrow \mathbb{R}$, we will denote ${w}^*$ its upper semi-continuous regularization, defined in $A$ by $\displaystyle w^*(x):=\limsup_{y\rightarrow x, y \in A}w(y)$.\\
We deduce from Lemma \ref{local} the following result :

\begin{lem}\label{sup}Let $T$ be a closed positive current defined in an open set $U\subset\mathbb{C}^n$. Let $(u_n)_n$ a locally uniformly bounded from above sequence of continuous co-$T$-psh functions. Then the function $\displaystyle {\sup_n}^*u_n$ is $T$-psh.\\
Moreover, given an open set $W\subset\mathbb{C}^n$ and a holomorphic map $f:W\rightarrow U$, we have $$\displaystyle ({\sup_n}^*u_n)\circ f = {\sup_n}^*(u_n\circ f).$$ \end{lem}

\noindent \textit{Proof} Fix an open set $V\subset\subset U$. For each $n\in\mathbb{N}$, there exists a sequence $(v_{n,i})_i$ of continuous plurisubharmonic functions defined in $V$, decreasing towards $u_n$ on $\mathcal{S}\cap V$. Thanks to Dini's Lemma, we may suppose $\displaystyle |v_{n,i}-u_n|<\frac{1}{n}$ on $\mathcal{S}\cap\overline{V}$ for all $n,i\in\mathbb{N}$. Then the sequence of plurisubharmonic functions $\displaystyle ({\sup_n}^* v_{n,i})_i$ decreases on $\mathcal{S}\cap V$ (the upper regularization is done here in $V$).
We will prove that its limit on $S\cap V$ is precisely $\displaystyle {\sup_n}^*({u_n}{|_{\mathcal{S}\cap V}})$, which suffices to conclude. 

For this purpose, it suffices to prove that $\forall i\in\mathbb{N}$, $\displaystyle ({\sup_n}^* v_{n,i}){|_{\mathcal{S}\cap V}}= {\sup_n}^*( v_{n,i|{\mathcal{S}\cap V}})$. Suppose, by contradiction, that there exists $x_0\in \mathcal{S}\cap V$ and $i\in\mathbb{N}$ such that $\displaystyle ({\sup_n}^* v_{n,i}){|_{\mathcal{S}\cap V}}(x_0)> {\sup_n}^* (v_{n,i}|_{\mathcal{S}\cap V})(x_0)$. Then by Proposition 4.8.2 in \cite{K}, there exists a plurisubharmonic function $v:V\rightarrow\mathbb{R}$ such that $\displaystyle\limsup_{x\rightarrow x_0,\text{ }x\in \mathcal{S}\cap V\setminus \{x_0\}}v(x)<v(x_0)$. But this contradicts Corollary \ref{maximum} applied to the $T$-psh function $v{|_{\mathcal{S}\cap V}}:\mathcal{S}\cap V\rightarrow \mathbb{R}$.

An adaptation of the idea of the proof of Proposition 2.9.16 in \cite{K} establishes the second statement, as follows. Denote by $J_f$ the Jacobian of $f$. For all $i\in\mathbb{N}$, the functions $\displaystyle ({\sup_n}^*u_{n,i})\circ f$ and $\displaystyle {\sup_n}^*(u_{n,i}\circ f)$ (the upper regularization is done here in $W$) are two plurisubharmonic functions defined in $W$ which are equal in $W\setminus \{J_f=0\}$, hence everywhere. The previous argument shows that they decrease on $\mathcal{S}\cap W$ when $i\rightarrow +\infty$, respectively towards the functions $\displaystyle ({\sup_n}^*u_n)\circ f $ and $\displaystyle {\sup_n}^*(u_n\circ f)$, which are therefore equal.$\blacksquare$\\ 

\noindent The following lemma gives an alternative definition for continuous $T$-psh functions.

\begin{lem}\label{lem.jensen} Let $T$ be a closed positive $(1,1)$-current with support $\mathcal{S}$ in an open set $U\subset\mathbb{C}^n$, and let $u:\mathcal{S}\rightarrow \mathbb{R}$ be a continuous function. Then $u$ is $T$-psh if and only if $\displaystyle u(x)\leq\int_{U}u d\Lambda_x$, $\forall x\in \mathcal{S}$ and $\forall\Lambda_x\in PJ_{\mathcal{S}}(x)$. 
\end{lem}

\noindent \textit{Proof} \noindent{\bf Necessary condition.} Consider a sequence $(v_n)_n$ of plurisubharmonic functions defined in an open set containing $\mathcal{S}$, such that $(v_n)_n$ decreases towards $u$ on $\mathcal{S}$. Let $x\in \mathcal{S}$ and $\Lambda_x\in PJ_{\mathcal{S}}(x)$. The monotone convergence theorem gives \begin{equation}\label{gnagnagna}\displaystyle v(x)=\lim_{n\rightarrow +\infty}v_n(x)\leq \lim_{n\rightarrow +\infty}\int_{U} v_n d\Lambda_x=\int_{U} v d\Lambda_x.\end{equation}

\noindent{\bf Sufficient condition.} Let $V\subset\subset U$ be an open set intersecting $\mathcal{S}$, and let $K:=\mathcal{S}\cap\overline{V}$. First let us recall an abstract duality theorem due to Edwards \cite{E}. Let $Z$ be a compact metric space, and $\mathcal{R}$ a cone of continuous functions from $Z$ to $\mathbb{R}$, such that $\mathcal{R}$ separates points and contains the constants. Given $x\in Z$, we denote by $I_x$ the set of probability measures $\Lambda$ on $Z$ which satisfy $\displaystyle v(x)\leq\int_{K}v d\Lambda$ for all $v\in \mathcal{R}$. Let $\phi:Z\rightarrow \mathbb{R}$ be a continuous function. Then Edwards' Theorem ensures that for all $ x\in Z$,
$$
\displaystyle\sup\left\{ v(x):v\in \mathcal{R}, v\leq \phi\right\} =\inf\left\{ \int_Z \phi d\Lambda : \Lambda\in I_x\right\}.
$$

\noindent Now let $u:\mathcal{S}\rightarrow \mathbb{R}$ be a continuous function such that $\forall x\in \mathcal{S}$, $\forall\Lambda_x\in PJ_{\mathcal{S}}(x)$, $\displaystyle u(x)\leq\int_{U}u d\Lambda_x$. Let also $(A_n)_n$ be a compact exhaustion of $U$, and fix $m\in\mathbb{N}$. Denote $\displaystyle K_{m}:=A_m\cap\mathcal{S}$. We apply Edwards' theorem with $\mathcal{R}$ being the set of continuous $K_{m}$-psh functions, $Z=K_{m}$ and $\phi=u$. Note that $I_x=PJ_{K_{m}}(x)$. Indeed, (\ref{gnagnagna}) implies that $PJ_{K_{m}}(x)\subset I_x$, while the inverse inclusion is trivial. We have $\forall x\in K_{m}$, $$
\displaystyle \sup\left\{ v(x)|v\in \mathcal{R}, v\leq u\right\} \leq u(x)\leq \inf\left\{ \int_{K_m} u d\Lambda |\Lambda\in I_x\right\}.
$$
Then, Edwards' theorem gives $\forall x\in K_{m}$, $$\displaystyle u(x)=\sup\left\{ v(x)|v\in R, v\leq u\right\}=\inf_{\Lambda_x\in PJ_{K_{m}}(x)}\int_{K_m}u d\Lambda_x.$$
We deduce that $\displaystyle u=u^* = {\sup}^*\left\{ v|v\in \mathcal{R}, v\leq u\right\}$ on $K_{m}$. But thanks to Lemma 2.3.4 in \cite{K} and our Lemma \ref{sup}, $\displaystyle {\sup}^*\left\{ v|v\in \mathcal{R}, v\leq u\right\}$ is $T$-psh in $\displaystyle \mathring{K}_{m}$. The result follows by taking the limit $m\rightarrow +\infty$.
$\blacksquare$\\ 

\subsection{Expression for relative Green's functions and applications}\label{subsectiontwo}
 
Denote by $\Delta\subset \mathbb{C}$ the unit disk of $\mathbb{C}$. Let $U\subset\mathbb{C}^n$ be an open set. Recall that an \textit{analytic disk in $U$ with center $x\in U$} is a continuous function $h:\overline\Delta\rightarrow \mathbb{C}^n$, analytic in $\Delta$, such that $h(0)=x$. If $\lambda$ is the normalized Lebesgue measure on $\partial \Delta$, an \textit{analytic-disk measure on $U$ for x} is a measure of the form $g_* \lambda$, where $g:\overline\Delta\rightarrow U$ is an analytic disk in $U$ with center $x$.

We will denote by $\Lambda_U(x)$ the set of analytic disks which are pluri-Jensen measures relative to $x$ with support in an open set $U$. By Theorem 6.3 in \cite{Ra}, $\displaystyle PJ_U(x)=\overline{\Lambda_U(x)}$.

We will need the following theorem of Poletsky (Corollary p. 170 in \cite{Po} , see also \cite{Po2}), generalized by Rosay (\cite{Ro}). Let $U$ be a connected complex manifold of dimension $\geq 1$. We denote by $\mathcal{H}_{z,U}$ the set of holomorphic functions $h:V_h\rightarrow U$ from some neighborhood $V_h$ of $\overline{\Delta}=\{z\in \mathbb{C}, |z|\leq 1\}$ (possibly depending on $h$) into $U$ such that $h(0)=z$. We also denote by $PSH(U)$ the set of plurisubharmonic functions defined on $U$. 

\begin{thm}\label{poletsky} Let $u:U\rightarrow \mathbb{R}$ be an upper semi-continuous function. With the previous notations, the function defined in $U$ by $$\displaystyle \tilde{u}(z):=\frac{1}{2\pi}\inf_{f\in \mathcal{H}_{z,U}}\int_0^{2\pi} u(f(e^{i\theta}))d\theta ,$$ if it is not everywhere equal to $-\infty$, belongs to $PSH(U)$ and verifies $\tilde{u}\leq {u}$. Moreover, this function $\tilde{u}$ is maximal among all the functions in $PSH(U)$ verifying this inequality.\\
If $u$ is only assumed borelian, we still have $\tilde{u}\leq u$, and $\tilde{u}^*$ is plurisubharmonic.\end{thm}

The last statement in the above theorem is the Lemma 2.2.12 in \cite{St}.
We prove below a version of this theorem, useful in our context. Let us recall some useful definitions.\\
An open bounded set $U\subset\mathbb{C}^n$ is called {\it{hyperconvex}} if there exists a continuous plurisubharmonic function $\rho:U\rightarrow ]-\infty,0[$ such that $\forall c\in ]-\infty,0[$, $\{z \in \mathbb{C}^n: \rho(z)<c\}\subset \subset U$. In order to simplify the presentation, we will also suppose here that a hyperconvex set is connected.\\
Let $U\subset\mathbb{C}^n$, $n\geq 1$, be a bounded open set, and let $A\subset\mathbb{C}^n$ be a compact set. The {\it{Green's function of $A\subset U$ relative to $U$}} is defined by $$\displaystyle G_{A,U}:= {\sup}^* \big\{ v \in PSH(U):\text{ } v\leq 0,\text{ } v|_{A}\leq -1\big\},$$where ${\sup}^*$ denotes the upper semi-continuous regularization of the upper envelope. A compact set $K\subset\mathbb{C}^n$ is called {\it{fat}} if $\displaystyle \overline{\mathring{K}}=K$. In this case, if $K$ is included into a bounded hyperconvex open set $U\subset\mathbb{C}^n$, the function $G_{K,U}$ is continuous by Proposition 4.5.3 of \cite{K}.\\

\begin{thm}\label{thetaclaque} Let $w:U\rightarrow \mathbb{R}$ be an upper semi-continuous function on a connected open set $U\subset \mathbb{C}^n$. Then the function defined in $U$ by \begin{equation}\label{eqone}\displaystyle\hat{w}(x):=\inf_{\mu_x\in PJ_U(x)}\int_U w d\mu_x\end{equation} is maximal among all the plurisubharmonic functions defined in $U$ and bounded from above by $w$ (if there are any).\\
In particular, if the open set $U$ is bounded and hyperconvex, and if $A\subset U$ is a fat compact set, the Green's function $A$ relative to $U$ equals 
\begin{equation}\label{eqtwo}\displaystyle G_{A,U}(x) =  \inf_{\mu_x\in PJ_U(x)}\int_U(-\mathbf{1}_A)^* d\mu_x=-\sup_{\mu_x\in PJ_U(x)}\mu_x(\mathring{A}).\end{equation}
In the case where the set $A\subset U$ is only assumed relatively compact, if the open set $U$ is bounded and hyperconvex, we have \begin{equation}\label{eqthree}\displaystyle G_{A,U}(x) \geq \inf_{\mu_x\in PJ_U(x)}\int_U(-\mathbf{1}_{\overline{A}}) d\mu_x.\end{equation}
\end{thm}

\noindent \textit{Proof} We denote again $\mathcal{H}_{z,U}$ the set of holomorphic functions $h:V_h\rightarrow U$ from a neighborhood $V_h$ of $\overline{\Delta}$ (possibly depending on $h$) into a neighborhood $U$ of $z\in\mathbb{C}^n$, such that $h(0)=z$. 

Let $\mu_x=g_*\lambda$ with $g:\overline{\Delta}\rightarrow \mathbb{C}^n$ an analytic disk. By Mergelyan's Theorem, we can find a sequence $(f_n)_n$ of $\mathcal{H}_{z,U}$ which converges to $g$ uniformly on $\overline{\Delta}$. 
The upper semi-continuity of $w$ implies that $\forall \theta\in[0,2\pi]$, $$\displaystyle\limsup_{n\rightarrow +\infty}w\circ f_n(\theta)\leq w\circ g(\theta).$$ 

\noindent It follows, by Fatou's Lemma, that : $$\displaystyle \limsup_{n\rightarrow +\infty} \frac{1}{2\pi}\int_0^{2\pi} w\circ f_n (e^{i\theta}) d\theta\leq \int_U w \mu_x.$$ 

\noindent (Note that this inequality also follows from Lemma 3.4.4 in \cite{K}). We deduce that $$\displaystyle\frac{1}{2\pi}\inf_{f\in \mathcal{H}_{z,U}}\int_{0}^{2\pi}w\circ f(e^{i\theta})d\theta\leq \inf_{\nu_x \in \Lambda_U(x)}\int_U w d\nu_x. $$ 

\noindent The inequality in the other direction being trivial, we finally have \begin{equation}\label{premiere}\inf_{\nu_x \in \Lambda_U(x)}\int_U w d\nu_x = \frac{1}{2\pi}\inf_{f\in \mathcal{H}_{z,U}}\int_{0}^{2\pi}w\circ f(e^{i\theta})d\theta.\end{equation}

\noindent On the other hand, recall that every pluri-Jensen measure is a limit of analytic disks. For $ \nu_x\in PJ_U(x)$, denote $(\nu_{x,n})_n$ a sequence of analytic disks, weakly convergent to $\nu_x$. Thanks to the upper semi-continuity of $w$ (see e.g. Lemma 2.8 in \cite{D}, or Lemma 3.4.4 in \cite{K}), we have $\forall \nu_x\in PJ_U(x)$, \begin{equation}\label{seconde}\displaystyle \inf_{\mu_x \in \Lambda_U(x)}\int_U w d\mu_x \leq\lim_{n\rightarrow +\infty}\int_U w d\nu_{x,n}\leq  \int_U w d\nu_x.\end{equation}
Bringing together Equations (\ref{premiere}) and (\ref{seconde}), we obtain $$\displaystyle \inf_{\nu_x \in PJ_U(x)}\int_U w d\nu_x \geq  \inf_{\nu_x \in \Lambda_U(x)}\int_U w d\nu_x = \frac{1}{2\pi}\inf_{f\in \mathcal{H}_{z,U}}\int_{0}^{2\pi}w\circ f(e^{i\theta})d\theta\geq \inf_{\nu_x \in PJ_U(x)}\int_U w d\nu_x,$$and hence the above inequalities are all equalities. Theorem \ref{poletsky} concludes the proof of (\ref{eqone}). 

Recall that $G_{A,U}$ is continuous by Proposition 4.5.3 of \cite{K}. Then (\ref{eqtwo}) comes from (\ref{eqone}), since $$\displaystyle G_{A,U}(x) = \inf_{\mu_x\in PJ_U(x)}\int_U(-\mathbf{1}_{\mathring{A}}) d\mu_x\geq \inf_{\mu_x\in PJ_U(x)}\int_U(G_{A,U}) d\mu_x=G_{A,U}(x).$$

\noindent {\it{Let us now show (\ref{eqthree}).}} Denote by $dist(\cdot,\cdot)$ the euclidian distance. Take $\epsilon>0$ such that $U_{\epsilon}:=\{z\in U, dist(z,\partial U)<\epsilon\}$ is nonempty and contains $A$. Take $x\in U_{\epsilon}$. Let $A_{m}:=\{z\in\mathbb{C}^n, dist(z,A)\leq \frac{1}{m}\}$ for every $m\in\mathbb{N}^*$ sufficiently large such that $A_m\subset\subset U_{\epsilon}$. Let then $(\nu_{x,n,k})_k$ be a weakly convergent sequence of measures in $ PJ_{{U_{\epsilon}}}(x)$ such that $$\displaystyle\inf_{\nu_x\in J_{{U_{\epsilon}}}(x)}\int_{{U_{\epsilon}}} \left(-\mathbf{1}_{A_n}\right)d\nu_x=\lim_{k\rightarrow +\infty }\int_{{U_{\epsilon}}} (-\mathbf{1}_{A_n})d\nu_{x,k,n}.$$
Denote $\nu'_{x,n}$ its limit, and note that $\nu'_{x,n}\in PJ_{\overline{U_{\epsilon}}}(x)$ (see Theorem 6.3 in \cite{Ra}).\\
Let also $\nu'_x\in PJ_{{\overline{U_{\epsilon}}}}(x)$ be an accumulation point of the sequence $(\nu'_{x,n})_n$, towards which we may suppose that it converges. Note that the sequence $(-\mathbf{1}_{A_{n}})_n$ is a non-decreasing sequence of lower semi-continuous functions. The Lemma 2.8 in \cite{D} used twice, implies that 
\begin{align*} 
 \displaystyle G_{A,U_{\epsilon}}(x) &\geq \lim_{n\rightarrow +\infty}G_{A_n,U_{\epsilon}}(z)=\lim_{n\rightarrow +\infty}\lim_{k\rightarrow +\infty}\int_{U_{\epsilon}}(-{\mathbf{1}}_{A_n})^* d\nu_{x,n,k} \\
 &\geq \limsup_{n\rightarrow +\infty}\limsup_{k\rightarrow +\infty}\int_{U_{\epsilon}}(-{\mathbf{1}}_{A_n}) d\nu_{x,n,k}
 \displaystyle\geq \limsup_{n\rightarrow +\infty}\int_{\overline{U_{\epsilon}}}(-\mathbf{1}_{A_n}) d\nu'_{x,n} \\
&\geq  
\int_{\overline{U_{\epsilon}}}(-\mathbf{1}_{\overline{A}}) d\nu'_{x} \geq \inf_{\mu_x\in J_{\overline{U_{\epsilon}}}(x)}\int_{\overline{U_{\epsilon}}}(-\mathbf{1}_{\overline{A}}) d\mu_x
\geq \inf_{\mu_x\in PJ_U(x)}\int_U(-\mathbf{1}_{\overline{A}}) d\mu_x.
\end{align*}
Now, since $U$ is assumed hyperconvex, Proposition 4.5.7 in \cite{K} allows us to pass to the limit $$\displaystyle G_{A,U}(x)=\lim_{\epsilon\rightarrow 0}G_{A,U_{\epsilon}}(x)\geq \inf_{\mu_x\in PJ_U(x)}\int_U(-\mathbf{1}_{\overline{A}}) d\mu_x.$$
$\blacksquare$\\

\begin{rmq}\label{rmq} It follows from Theorem \ref{thetaclaque} that for a bounded hyperconvex open set $U\subset \mathbb{C}^n$, if $A\subset U$ is a fat compact set such that $\hat{A}\subset U$, then $\forall x\in U$, $\displaystyle \sup_{\nu_x\in PJ_U(x)}\nu_x(\mathring{A}) =\sup_{\nu_x\in PJ_U(x)}\nu_x(\mathring{\hat{A}})$. Indeed, $\forall x\in U$, we have $$\displaystyle \inf_{\nu_x\in PJ_U(x)}\int_U(-\mathbf{1}_A)^* d\nu_x = G_{A,U}(x) = G_{\hat{A},U}(x)=\inf_{\nu_x\in PJ_U(x)}\int_U(-\mathbf{1}_{\hat{A}})^* d\nu_x.$$
(The equality $G_{A,U} = G_{\hat{A},U}$ comes from Theorem 5.1.7 in \cite{K} and Corollary 5.3.4 in \cite{K}.) Thus ${A}$ is polynomially convex if and only if $\forall x\notin {A}$, $$\displaystyle \sup_{\nu_x\in PJ_{U}(x)}\nu_x(\mathring{A})<1.$$\\
 Finally, since $PJ_{\overline{U}}(x)$ is closed (Lemma 7.2 in \cite{GR}, see also Theorem 6.3 in \cite{Ra}), Lemma 3.4.4 in \cite{K} implies $$\displaystyle A=\hat{A}\iff \forall x\notin {A},\text{ }\forall\nu_x\in PJ_U(x),\text{ }\nu_x({\mathring{A}})<1.$$\end{rmq} 

\begin{app} In particular, Remark \ref{rmq} allows to deduce that, if $u:\mathbb{C}^n\rightarrow \mathbb{R}^+$ is a continuous plurisubharmonic function, then the set $\{u = 0\}$ is polynomially convex, provided that it is bounded and fat.\\
Indeed, suppose that $\{u= 0\}$ is bounded and fat but not polynomially convex. By Remark \ref{rmq}, one can find $y\in\mathbb{C}^n\setminus \{u= 0\}$ and a measure $\mu_y\in PJ_y(\mathbb{C}^n)$ whose support is included in $\{u =0\}$. Then $$0< u(y)\leq \int_{\mathbb{C}^n} u d\mu_y\leq \int_{\mathbb{C}^n} 0 d\mu_y = 0.$$
We have thus proved by contradiction that $\{u = 0\}$ is polynomially convex, if it is bounded and fat. \end{app}

\begin{cor}\label{condpluripolar} Let $U\subset\mathbb{C}^n$ be an open set, and let $A\subset U$ be a complete pluripolar set, i.e. $A=\{u=-\infty\}$ for a plurisubharmonic function $u:U\rightarrow \mathbb{R}\cup \{-\infty\}$. Then $\forall x\in U\setminus A$, the set $A$ is of null measure for every pluri-Jensen measure in $PJ_{U}(x)$.\\
Conversely, suppose $U$ bounded and hyperconvex. If there exists $x\in U\setminus \overline{A}$ such that the set $A$ is of null measure for every pluri-Jensen measure in $PJ_{U}(x)$, then $A$ is pluripolar.\end{cor}

\noindent \textit{Proof} Let $ x\in\mathbb{C}^n\setminus A$, and $\mu_x\in PJ_{U}(x)$. We have $\mu_x(A)=0$, since $$-\infty<u(x)\leq \int_{U}u d\mu_x.$$
Conversely, we may suppose that $A$ is compact. Suppose $\exists x\in U\setminus \overline{A}$, $\forall \nu_x\in PJ_{U}(x)$, we have $\nu_x(A)=0$. Then Theorem \ref{poletsky} implies that $$\displaystyle 0=\inf_{\nu_x\in PJ_{U}(x)}\int_U\left(-\mathbf{1}_{\overline{A}}\right)d\nu_x\leq G_{A,U}(x).$$
The function $G_{A,U}$ is null in a point outside $\overline{A}$, hence everywhere on $U$, and thus $A$ is pluripolar. $\blacksquare$\\

\begin{rmq} Note that, given a borelian function $u:U\rightarrow \mathbb{R}\cup \{-\infty\}$ such that the set $\{u=-\infty\}$ is complete pluripolar, we deduce, by applying Theorem \ref{poletsky} to $U\setminus \{u=-\infty\}$, that the statement of  Theorem \ref{poletsky} remains valid when replacing $u:U\rightarrow \mathbb{R}$ by $u:U\rightarrow \mathbb{R}\cup \{-\infty\}$. Similarly, the statement of Theorem \ref{thetaclaque} remains valid when replacing $w:U\rightarrow \mathbb{R}$ by $w:U\rightarrow \mathbb{R}\cup \{-\infty\}$.\end{rmq}

\noindent We obtain the following lemma, which we will often use in the sequel :

\begin{lem}\label{restriction} Let $T$ be a closed positive $(1,1)$-current defined in a connected open set $U\subset\mathbb{C}^n$. Let $v:\mathcal{S}\rightarrow\mathbb{R}^-$ be a $T$-psh function. If $u:U\rightarrow \mathbb{R}^-$ is the extension of $v$ to $U$ with the value $0$, we have $\forall x\in \mathcal{S}$, $$\displaystyle \inf_{\nu_x\in PJ_{U}(x)}\int_U u d\mu_x =\inf_{\nu_x\in PJ_{\mathcal{S}}(x)}\int_{\mathcal{S}}v d\mu_x=v(x).$$ 
\end{lem}

\noindent \textit{Proof} First suppose that $v=w{|_{\mathcal{S}}}$ for some plurisubharmonic function $w:U\rightarrow \mathbb{R}$ such that $w{|_{\mathcal{S}}}<0$. Denote by $W$ the open set $\{w<0\}\subset U$, and note that $\mathcal{S}\subset W$. Then Theorem \ref{thetaclaque} and (\ref{eq.jensen}) imply that $\forall x\in\mathcal{S}$ we have
\begin{align*}
\displaystyle v(x)&= \inf_{\nu_x\in PJ_{U}(x)}\int_U w d\nu_x \leq\inf_{\nu_x\in PJ_{W}(x)}\int_U w d\nu_x\\
& \leq \inf_{\nu_x\in PJ_{W}(x)}\int_U{u}d\nu_x\leq \inf_{\nu_x\in PJ_{\mathcal{S}}(x)}\int_U{u}d\nu_x=v(x).
\end{align*}
Now suppose only that $v:\mathcal{S}\rightarrow \mathbb{R}^-$ is a $T$-psh function as in the statement. Let $(v_n)_n$ be a sequence of plurisuharmonic functions defined in $U$, decreasing on $\mathcal{S}$ towards $v$. Take $x\in \mathcal{S}$ and $\mu_x\in PJ_{U}(x)$. Take $\epsilon>0$. Then there exists $n_0$ such that $v_{n}(x)-\epsilon<0$ for $n\geq n_0$. There also exists an open neighborhood $O\subset U$ of $x$ such that $\displaystyle O\cap \mathcal{S}\subset \bigcap_{n\geq n_0}\{v_n-\epsilon<0\}$. Up to increasing $n_0$, we may suppose that the compact support of $\mu_x$ is included in $O$. Denote by $u_{n,\epsilon}$ the extension to $U$ with the value $0$ of ${(v_n-\epsilon)}{|_{\mathcal{S}}}$, and denote by $u_n$ the extension to $U$ with the value $0$ of ${v_n}{|_{\mathcal{S}}}$. By the previous step,$$\displaystyle v_n(x)-\epsilon=\inf_{\nu_x\in PJ_{O}(x)}\int_U u_{n,\epsilon}(x)\nu_x\leq \int_U u_n d\mu_x.$$ The monotone convergence theorem and the arbitrary choice of $\epsilon>0$ give $\displaystyle v(x)\leq\int_U u d\mu_x$. Then $\displaystyle v(x)=\inf_{\mu_x\in PJ_{U}(x)} \int_U u d\mu_x$ since $\delta_x\in PJ_U(x).$$\blacksquare$\\

\noindent Let us now introduce a new notion, which we will use in the next result. 
\begin{defn}
Let $\mathcal{S}$  be the support of a closed positive $(1,1)$-current $T$. A set $E\subset \mathcal{S}$ is called \emph{$T$-pluripolar} if there exists a $T$-psh function $v:\mathcal{S}\rightarrow \mathbb{R}\cup \{-\infty\}$ such that $e^u$ is continuous and $E\subset \{v=-\infty\}$.
\end{defn}

We prove the following result.
\begin{pro}\label{pluripolar} A set $F\subset S$ which is $T$-pluripolar is pluripolar.\end{pro}
\noindent \textit{Proof} Without loss of generality, we may suppose that $F$ is bounded, since a locally pluripolar set is pluripolar. Let $u:\mathcal{S}\rightarrow \mathbb{R}\cup\{-\infty\}$ be a $T$-psh function with $e^u$ continuous, and $E=\{u=-\infty\}\subset \mathcal{S}$ a $T$-pluripolar set containing $F$. Let $U\subset\mathbb{C}^n$ be an open bounded hyperconvex set containing $F$, and let $(u_n)_n$ be a sequence of continuous plurisubharmonic functions defined in $U$, decreasing towards $u$ on $\mathcal{S}$. Let $v_n:={u_{n}}{|_{\mathcal{S}}}$ extended to $U$ with the value $0$. Let $K\subset U\cap\mathcal{S}$ be a non-pluripolar compact set containing also $F$. Up to shrinking $U$, we may suppose that $u_n\leq 0$ on $U\cap \mathcal{S}$ for $n\geq n_0$. Theorem \ref{thetaclaque} and Lemma \ref{restriction} imply that $\forall n\geq n_0$, $\forall c\in\mathbb{R}^{-}\setminus \{0\}$, we have $\forall x\in U\cap\mathcal{S}$,
\begin{align*}& \displaystyle G_{\{u_n\leq c\}\cap K,U}(x)\geq \inf_{\mu_x\in PJ_U(x)}\int_{U}\left( -\mathbf{1}_{\{u_n\leq c\}\cap K}\right)d\mu_x\geq \inf_{\mu_x\in PJ_U(x)}\int_{\{u_n\leq c\}\cap\mathcal{S}}\left( -\mathbf{1}_{\{u_n\leq c\}}\right)d\mu_x\\
& \geq \frac{1}{|c|}\inf_{\mu_x\in PJ_U(x)}\int_{\{u_n\leq c\}\cap\mathcal{S}}v_n d\mu_x  \geq \frac{1}{|c|}\inf_{\mu_x\in PJ_U(x)}\int_{U\cap\mathcal{S}}v_n d\mu_x =\frac{u_n(x)}{|c|}.\end{align*}

\noindent Corollary 4.7.8 in \cite{K} allows us to pass to the limit :  $\forall c\in\mathbb{R}^{-}\setminus \{0\}$, $\forall x\in  U\cap\mathcal{S}$, $$\displaystyle G_{\{u\leq c\}\cap K,U}(y)=\lim_{n\rightarrow +\infty}G_{\{u_n\leq c\}\cap K,U}(y)\geq \frac{u(x)}{|c|}.$$Now let $y\in (U\cap\mathcal{S})\setminus E$. We have $$\displaystyle G_{E\cap K,U}(y)\geq\lim_{c\rightarrow -\infty}G_{\{u\leq c\}\cap K,U}(y)\geq 0.$$
Since we obtain $G_{E\cap K,U}=0$ on $(U\cap\mathcal{S})\setminus E$, hence on $U$, it follows that the set $E\cap K$ is pluripolar. Thus $F\subset E\cap K$ is pluripolar.$\blacksquare$\\

\begin{rmq}\label{remarque} It is not always true that a pluripolar set is $T$-pluripolar : for instance, $B(0,1)\cap\{z_1 = 0\}$ is pluripolar in $\mathbb{C}^2 = \{(z_1,z_2): z_1, z_2 \in \mathbb{C}\}$. This set is not $T$-pluripolar for $T = dd^c \log{ |z_1|}$, because $(z_1,z_2) \mapsto \log |z_1|\notin L^1_{loc}(\|T\|)$. \\
Using Proposition \ref{pluripolar}, the $T$-pluripolar sets are the pluripolar subsets of $supp(T)$ which are included in a set of the form $\{u=-\infty\}$, for a continuous plurisubharmonic function $u\in L^1_{loc}(\|T\|)$. This is the case, for example, if the set $\{u=-\infty\}\cap supp(T)$ is relatively compact in a ball contained in the domain of definition of $u$, by Proposition 3.1 in \cite{D}, or if it is of Hausdorff measure $\mathcal{H}_1(E)=0$, by Theorem 3.5 in \cite{D}. \end{rmq}

\noindent We will further need the following lemma :

\begin{lem}\label{proximite}
\noindent Let $U\subset \mathbb{C}^n$ be an open bounded hyperconvex set, and $A\subset\subset U$ a non pluripolar fat compact set. Then, for every $\epsilon>0$, there exists a neighborhood $V\subset\subset U$ of $A$ such that for every $x\in V$, there exists a pluri-Jensen measure $\nu_x\in PJ_{ U}(x)$ with $\displaystyle\nu_x\left(A\right)\geq 1-\epsilon$.\end{lem}

\noindent \textit{Proof} By Theorem \ref{thetaclaque}, $$\displaystyle\inf_{\nu_x\in PJ_U(x)}\int_U(-\mathbf{1}_A) d\nu_x\leq G_{A,U}(x).$$ It suffices thus to take $V:=\{G_{A,U}<-1+\frac{\epsilon}{2}\}$.$\blacksquare$\\ 

\subsection{Properties analogous to the ones for plurisubharmonic functions}\label{similar}

We prove here several results which transpose to $T$-psh functions certain properties of plurisubharmonic functions.

\begin{lem}\label{trucdeouf} Let $T$ be a closed positive current in an open subset of $\mathbb{C}^n$. Let $(u_m)_m$ be a locally uniformly bounded from above sequence of continuous co-$T$-psh functions. Suppose that the upper semi-continuous regularization $u^*_{\infty}:\mathcal{S}\rightarrow \mathbb{R}$ of the function $u_{\infty}:=\displaystyle \limsup_{m\rightarrow +\infty}u_m$ is continuous. Then the set $\{u_{\infty}<u^*_{\infty}\}$ is pluripolar.\end{lem}

\noindent \textit{Proof} {\it{First suppose that $u^*_{\infty}$ is constant.}}  Up to replacing $u_m$ by $u_m-u^*_{\infty}$, we may suppose that $u^*_{\infty}=0$. Let $V$ be an open hyperconvex subset of $\mathbb{C}^n$ intersecting $\mathcal{S}$, and $K\subset V\cap\mathcal{S}$ a non-pluripolar compact set.  Let a constant $\epsilon>0$. For every $m\in\mathbb{N}$, by Dini's Lemma, up to shrinking $V$, there exists a continuous plurisubharmonic function $v_m:V\rightarrow\mathbb{R}$, with $v_m\geq u_m$ on $\mathcal{S}$ and $|v_m-u_m|\leq \frac{1}{m}$ on $\mathcal{S}$. We can also suppose that for every $m\in\mathbb{N}$ we have $v_m<\epsilon$ on $\mathcal{S}\cap V$. Then $\displaystyle u_{\infty}= \limsup_{m\rightarrow +\infty}v_{m}|_{\mathcal{S}}$ on $V\cap\mathcal{S}$.

\noindent Let a constant $c<-\epsilon$. Let $x\in V\cap\mathcal{S}\setminus \{u_{\infty}<0\}$. There exists an extracted subsequence $(v_{m_j})_j$ of $(v_m)_m$ such that $(v_{m_j}(x))_j$ converges to $0$. Then, $\forall j\in \mathbb{N}$, a computation similar to the one in the proof of Proposition \ref{pluripolar} gives $\forall j\in\mathbb{N}$, $$\displaystyle G_{\{\sup_{m\geq m_j} v_{m}-\epsilon\leq c\}\cap K,V}(x)\geq G_{\{v_{m_j}-\epsilon\leq c\}\cap K,V}(x)\geq \frac{v_{m_j}(x)-\epsilon}{|c|}.$$ Corollary 4.7.8 in \cite{K} allows us to pass to the limit $$\displaystyle G_{{\{u_{\infty}-\epsilon\leq c\}}\cap K,V}(x)= \lim_{j\rightarrow +\infty}G_{\{\sup_{m\geq m_j} v_{m}-\epsilon\leq c\}\cap K,V}(x)    \geq  \lim_{j\rightarrow +\infty}\frac{v_{m_j}(x)-\epsilon}{|c|}=\frac{-\epsilon}{|c|}.$$
Corollary 4.7.8 in \cite{K} allows us to take the limit $\epsilon\rightarrow 0$, and to obtain \begin{equation}\label{estpluripolaire}G_{{\{u_{\infty}\leq c\}\cap K},V}(x)=0.\end{equation}
The function $G_{{\{u_{\infty}\leq c\}\cap K},V}$ is null at $x$, hence on $V$. Finally, Corollary 4.7.8 in \cite{K} allows again to pass to the limit $$\displaystyle G_{{\{u_{\infty}< 0\}\cap K},V}= \lim_{c\rightarrow 0^-}G_{{\{u_{\infty}\leq c\}\cap K},V}=0.$$
Thus, by the arbitrary choice of $K$, the set $\{u_{\infty}<0\}$ is pluripolar.

\noindent {\it{Now we suppose only that $u^*_{\infty}$ is continuous.}} Remark that in this case the previous reasoning shows that for every compact set $K\subset\mathcal{S}$ and every real number $c<u^*_{\infty}$ on $K$, the set $\{u_{\infty}<c\}\cap K$ is pluripolar. Indeed, under these conditions, (\ref{estpluripolaire}) is still true. \\
Let $(s_m)_m$ be an increasing sequence of step functions, locally uniformly convergent towards $u^*_{\infty}$. By the previous remark, for all $m\in\mathbb{N}$, the set $\displaystyle\{\phi_{\infty}<s_m\}$ is pluripolar. Thus the set $\displaystyle\{\phi_{\infty}<\phi_{\infty}^*\}=\bigcup_{m\in\mathbb{N}}\{\phi_{\infty}<s_m\}$ is pluripolar.$\blacksquare$\\

We further present a version of Hartogs's Theorem for $T$-psh functions. 

\begin{thm}\label{hartogs} Let $\mathcal{S}$ be the support of a closed positive $(1,1)$-current $T$ in an open bounded set $U\subset \mathbb{C}^n$. Let also $Q\subset\mathbb{C}^n$ be a fat compact set such that $K:=Q\cap \mathcal{S}$ is non-pluripolar, and $E\subset K$ a set whose closure $\overline{E}$ in $K$ is of empty interior for the induced topology.\\
Let $(u_n)_n:U\cap \mathcal{S}\rightarrow \mathbb{R}$ be a sequence of continuous co-$T$-psh functions uniformly bounded from above on $U$, such that $\displaystyle \lim_{n\rightarrow +\infty}\left(\sup_{i\geq n}u_i\right)^*\leq f$ on $K\setminus E$ for some continuous function $f:\mathcal{S}\cap U\rightarrow \mathbb{R}$. Then for every $\epsilon>0$, there exists $n_0\in\mathbb{N}$, and a neighborhood $V\subset U$ of $K$ such that $u_n\leq f+\epsilon$ on $V\cap \mathcal{S}$ for all $n\geq n_0$. In particular, $\displaystyle \limsup_{n\rightarrow +\infty}u_n\leq C$ on $K$. \end{thm}

\noindent \textit{Proof} {\it{First suppose that $f$ is equal to a constant $C$.}} We may suppose, up to subtracting a constant, that $u_n< 0$ in $U\cap\mathcal{S}$, and that $C<0$. Let us denote $\displaystyle u'_n:=(\sup_{i\geq n}u_i)^*$. Then, the sequence of upper semi-continuous functions $(u'_n)_n$ decreases on $\mathcal{S}$ towards $\displaystyle u:=\lim_{n\rightarrow +\infty}\left(\sup_{i\geq n}u_i\right)^*$.\\
By Lemma \ref{proximite}, there exists a neighborhood $V$ of $K$ such that $\forall x\in V\cap\mathcal{S}$, we can find $\nu_x\in PJ_{U}(x)$ such that $\nu_x(K)>1-\frac{\epsilon}{2|C+\epsilon|}$.

By Egoroff's Theorem, there exists a set $A\subset K$, with $\nu_x(A)<\frac{\epsilon}{2|C+\epsilon|}$, such that $(u'_n)_n$ converges uniformly towards $u$ on $K\setminus A$. Fix $\epsilon>0$ such that $C+2\epsilon<0$. Then $\exists n_0\in\mathbb{N}$ such that $\forall z\in K\setminus (A\cup E)$, $\forall n\geq n_{0}$ : $\displaystyle u_n(z)\leq {u'_n}(z)\leq  C+\epsilon.$ The continuity of the functions $u_n$ then implies that $\forall z\in K\setminus A,\text{ }\forall n\geq n_{0}$ we have $$u_n(z)\leq C+\epsilon.$$ Extend each $u_n$ to $U$ with the value $0$. By the negativity of the functions $u_n$ on $\mathcal{S}\cap U$, we have for $n\geq n_{0}$ and $x\in V\cap\mathcal{S}$, thanks to (\ref{eq.jensen}) and Lemma \ref{restriction}, \begin{align*}
 \displaystyle u_n(x)&=\inf_{\mu_x\in PJ_{U\cap\mathcal{S}}(x)}\int_U u_n d\mu_x=\inf_{\mu_x\in PJ_{U}(x)}\int_U u_n d\mu_x\leq \int_U u_n d\nu_x\\
  &\leq  \int_{K\setminus A} u_n d\nu_x\leq \nu_x(K\setminus A)\sup_{K\setminus A} u_n\leq \left(1-\frac{\epsilon}{|C+\epsilon|}\right)(C+\epsilon)=C+2\epsilon. 
  \end{align*}
  
\noindent {\it{Now suppose that $f$ is a continuous function.}} Let $(s_m)_m$ be a sequence of step functions decreasing towards $f$, uniformly on $K$. There exists $m_0\in \mathbb{N}$ such that $s_{m_0}-f\leq \epsilon$. The previous step implies that, for $m$ sufficiently large we have on $K$ : $$u_m\leq s_{m_0}\leq f+\epsilon,$$
which proves the stated result.
$\blacksquare$\\

We deduce the following result :

\begin{pro}\label{uegalec} Let $\mathcal{S}$ be the support, supposed non-pluripolar, of a closed positive $(1,1)$-current $T$ in an open set $U\subset \mathbb{C}^n$. A $T$-psh function $u:\mathcal{S}\rightarrow \mathbb{R}$, constant outside a set $E$ whose closure in $\mathcal{S}$ is of empty interior, is constant everywhere.\end{pro}

\noindent \textit{Proof} Let us suppose that $u=C$ outside $E$. The upper semi-continuity of $u$ implies that $E=\{u\neq C\}=\{u>C\}$. Let $(u_n)_n$ be a sequence of continuous plurisubharmonic functions decreasing on $\mathcal{S}$ towards $u$. Theorem \ref{hartogs} applied to $(u_n)_n$ on compact subsets of $\mathcal{S}$ concludes.$\blacksquare$\\ 

This result allows us to weaken the conditions required in the statement of  Theorem \ref{hartogs} :

\begin{cor}\label{hartogs2} In Theorem \ref{hartogs}, if the current $T$ does not charge pluripolar sets, we can replace the condition $\displaystyle \lim_{n\rightarrow +\infty}\left(\sup_{i\geq n}u_i\right)^*\leq f$ by the less restrictive condition $\displaystyle \lim_{n\rightarrow +\infty}\left(\sup_{i\geq n}u_i\right)\leq f$.\end{cor}

\noindent \textit{Proof} We use the same notations as those in the statement of Theorem \ref{hartogs}. As in the proof of this theorem, we can suppose that the continuous function $f$ is equal to a constant $C\in\mathbb{R}$. It suffices to show that $$\displaystyle \lim_{n\rightarrow +\infty}\sup_{i\geq n}u_i\leq C\Rightarrow\lim_{n\rightarrow +\infty}\left(\sup_{i\geq n}u_i\right)^*\leq C.$$

It is sufficient again, by Lemma \ref{sup} and Proposition \ref{uegalec}, to show that $\displaystyle \lim_{n\rightarrow +\infty}\sup_{i\geq n}u_i\leq C$ implies that $\displaystyle\lim_{n\rightarrow +\infty}\max\left(\left(\sup_{i\geq n}u_i\right)^*, C\right)=C$
on $\mathcal{S}\setminus E$, where $E\subset \mathcal{S}$ is a closed set of empty interior for the induced topology. Suppose that the last condition is false. By Baire's theorem, there exists $k\in\mathbb{N}\setminus \{0\}$ such that $$\displaystyle\lim_{n\rightarrow +\infty}\left(\sup_{i\geq n}u_i\right)^*\geq C+\frac{1}{k}$$ in an open set (for the induced topology) $W\subset\mathcal{S}$. Then, $\displaystyle\left(\sup_{i\geq m}u_i\right)^*\geq C+\frac{1}{k}$ in $W$, $\forall m\in\mathbb{N}$. Hence there exists a nonempty set $O\subset W$ dense in $W$ (for the induced topology) such that $\forall m\in\mathbb{N}$, $\displaystyle \sup_{i\geq m}u_i \geq C+\frac{1}{2k}$ on $O$. This implies that $\displaystyle \lim_{n\rightarrow +\infty}\sup_{i\geq n}u_n> C$ on $O$, which finishes the proof.$\blacksquare$\\ 

We can now state a result which improves Corollary \ref{maximum} for certain currents :

\begin{pro}\label{liouville} Let $T$ be a closed positive $(1,1)$-current in an open set $U\subset\mathbb{C}^n$, with support $\mathcal{S}\subset U$. Let $(v_m)_m$ be a sequence of $T$-psh functions, such that $\displaystyle\limsup_{m\rightarrow +\infty} v_m\leq C$ for some $C\in\mathbb{R}$. If $\displaystyle\lim_{m\rightarrow +\infty}v_m(x)=C$ for some $x\in\mathcal{S}$, then $(v_m)_m$ converges towards $C$ outside a pluripolar set in $\mathcal{S}$.\\ In particular, if $\|T\|$ does not charge pluripolar sets, then any $T$-psh function that reaches its maximum is constant.
\end{pro}

\noindent \textit{Proof} Let us suppose $C=0$. Note that the function $\displaystyle {\limsup_{m\rightarrow +\infty}}^*{v_m}:\mathcal{S}\rightarrow \mathbb{R}$ is $T$-psh by Lemma \ref{sup}. Extend each function $v_m$ to $U$ with the value $0$. Fatou's lemma and Lemma \ref{restriction} imply that $\forall\nu_x\in PJ_{U}(x)$, $$\displaystyle 0=\lim_{m\rightarrow +\infty}v_m(x)\leq{\liminf_{m\rightarrow +\infty}}\int_U{v_m}  d\nu_x \leq{\limsup_{m\rightarrow +\infty}}\int_U{v_m}  d\nu_x \leq\int_U{\limsup_{m\rightarrow +\infty}}^*{v_m} d\nu_x \leq 0.$$ Finally $\forall\nu_x\in  PJ_{U}(x)$, $\displaystyle\lim_{m\rightarrow +\infty}\int_U{v_m}  d\nu_x=0$. Thus $(v_m)_m$ extended to $U$ with the value $0$ converges $\nu_x$-a.e. towards $0$ on the support of every measure $\nu_x\in PJ_{U}(x)$. Since $PJ_{U}(x)$ is invariant by rotations and homotheties with center $x$, Lemma \ref{local} and Corollary \ref{condpluripolar} imply that $(v_n)_n$ converges towards $0$ on $\mathcal{S}\setminus E$, where $E\subset\mathcal{S}$ is a pluripolar set. In particular, if $\|T\|$ does not charge pluripolar sets, then any $T$-psh function that reaches its maximum is constant outside such a set $E$, hence everywhere by Proposition \ref{uegalec}.$\blacksquare$\\

\noindent Finally, let us prove an inequality of Chern-Levine-Nirenberg type for $T$-psh functions :

\begin{pro}\label{tchern} Let $r_0>0$ be a real number. If $u$ is a continuous $T$-psh function, then there exists a constant $c\geq 0$ such that $\forall r< \frac{r_0}{2}$, $$\displaystyle T\wedge dd^c u (B(z,r))\leq \frac{c}{r^{2n}}\|T\|\left(B(z,2r)\right)\sup_{supp(T)\cap B(z,2r)}u.$$\end{pro} 

\noindent \textit{Proof} It suffices to prove the above formula when $u$ is the restriction of a plurisubharmonic function to the support $\mathcal{S}$ of $T$. The monotone convergence theorem then allows us to extend the formula to the general case.\\
Let $\xi$ be a positive step function equal to $1$ in $B(z,1)$ and to $0$ in $B(z,2)$. Let us denote $\xi_r(z):=\xi(\frac{z}{r})$. Then $\forall r\in ]0,1[$,  $$\displaystyle  T\wedge dd^c u (B(z,r))\leq \int_{\mathbb{C}^n}\xi_r T\wedge dd^c u = \int_{\mathbb{C}^n}u T\wedge dd^c \xi_r \leq \frac{c}{r^{2n}}\|T\|\left(B(z,2r)\right)\sup_{\mathcal{S}\cap B(z,2r)} u.$$
$\blacksquare$\\   

\section{Applications to dynamics in $\mathbb{C}^n$}\label{dynamics}

\subsection{Background}

Given a polynomial application $f:\mathbb{C}^n\rightarrow \mathbb{C}^n$, we extend it to a meromorphic endomorphism $[P_1,\ldots,P_n] : \mathbb{P}^n\rightarrow\mathbb{P}^n$, where $P_1,\ldots,P_n$ are homogeneous  polynomials. The choice of $\mathbb{P}^n$ as compactification of $\mathbb{C}^n$ is practical in our case because the group $H^{1,1}(\mathbb{P}^n,\mathbb{R})$ is of dimension $1$, and the classes which compose it are K\"ahlerian, in the sense that they contain a strictly positive smooth $(1,1)$-form. Denote by $H\subset\mathbb{P}^n$ the hyperplane at infinity. Denote by $I_f\subset H$ the {\it{indetermination set}} of $f$, i.e. the set of points in the neighborhood of which $f$ is not holomorphic. We denote $\displaystyle X_f:=\overline{f(H\setminus I_f)}\subset H$. This is an analytic set (see for example \cite{G} for details). 

Let $\omega$ be a K\"ahlerian form of mass 1 on $\mathbb{P}^n$. We define the \textit{dynamical degrees} $\lambda_1(f),\ldots,\lambda_n(f)\geq 1$ of $f$ by $$\lambda_k(f):=\liminf_{m\rightarrow +\infty}\bigg(\int_{\mathbb{C}^n}(f^{m*}\omega^k)\wedge\omega^{n-k}\bigg)^{\frac{1}{m}}.$$
Note that they are independent of the choice of the K\"ahler form $\omega$. When there is no risk of ambiguity, we will simply denote them by $\lambda_k$, and put $\lambda:=\lambda_1$.

The pull-back operation of a current induces a linear action $f^*$ on the cohomology spaces, defined by $f^*\{ S\}:=\{f^* S\}$. A polynomial application $f:\mathbb{C}^n\rightarrow \mathbb{C}^n$ extended to $\mathbb{P}^n$ is called {\it{algebraically stable}} if the linear action induced by $f^*$ on the cohomology spaces $H^{k,k}(\mathbb{P}^n,\mathbb{R})$ (of dimension 1) commutes with the composition operation, i.e. when $(f^n)^* = (f^*)^n$ for all $n\in \mathbb{N}$. When $f$ is algebraically stable, the spectral radius of the linear action induced in cohomology $f^* : H^{k,k}(\mathbb{P}^n,\mathbb{R})\rightarrow H^{k,k}(\mathbb{P}^n,\mathbb{R})$ equals the $k$-th dynamical degree $\lambda_k$. By duality, it also equals the spectral radius $\lambda'_{n-k}$ of the action induced by $f_*$ on $H^{n-k,n-k}(\mathbb{P}^n,\mathbb{R})$. In order for $f$ to be algebraically stable, it suffices to have $I_f\cap X_f=\emptyset$.

We will be interested in polynomial applications verifying $\lambda_1>\lambda_2$, whose properties will be seen to generalize those of the H\'enon's applications of $\mathbb{C}^2$. The concavity inequalities fulfilled by the $\lambda_i$ (Theorem 2.4.a in \cite{G}) imply that $\lambda_1$ is a dynamical degree strictly dominant ($\displaystyle\lambda_1>\max_{j\neq 1}\lambda_j$).\\

Given an algebraically stable polynomial endomorphism $f$ of $\mathbb{P}^n$ of dynamical degree $\lambda_1$ strictly dominant, and $\omega$ a K\"ahlerian form of mass $1$, the sequence  $\displaystyle\frac{1}{{\lambda_1}^m}f^{m*}\omega$ converges, in the weak sense of currents, towards a closed positive $(1,1)$-current $T^+$ (see e.g. Theorem 1.1 in \cite{B} for a general statement). The current $T^+$ is invariant by $\frac{1}{\lambda_1}f^*$ and does not depend on the choice of $\omega$. We denote by $J$ its support. 

If moreover $f$ is an automorphism, we can define a canonical positive $(n-1,n-1)$-current $\frac{1}{\lambda_1}f_*$-invariant $\displaystyle T^-:= \lim_{m\rightarrow +\infty}  \frac{1}{{\lambda_1}^m}{f^m}_*\omega^{n-1}$. The measures $\|T^+\|$ and $\|T^-\|$ do not charge pluripolar sets.

More generally, even if $f$ is not an automorphism, the current $T^-$ is still well defined if $f$ satisfies $\displaystyle \lim_{z\in \mathcal{C}_f,\text{ }\|z\|\rightarrow +\infty}f(z)\in X_f$, where $\mathcal{C}_f$ is the critical set of $f$ (see section 2.2 in \cite{G2}, see also \cite{GS2}).

The canonical invariant measure $T^+\wedge T^-$ is then well defined. It does not charge pluripolar sets. We will denote by $J^-$ the support of such current $T^-$. 

If $f$ is a polynomial endomorphism of $\mathbb{P}^n$ which is algebraically stable and verifies $\lambda_1>\lambda_2$, then there exists a plurisubharmonic function $G$ defined on $\mathbb{C}^n\subset\mathbb{P}^n$, called {\it{dynamical Green's function of $f$}}, such that $T^+=dd^c G$. It can be written as the limit in $L^1_{loc}$ of the sequence of term $ \frac{1}{{\lambda_1}^m}log(1+\|f^m\|)$. The function $G$ is continuous and pluriharmonic outside the support $J$ of $T^+$.

The current $T^+$ intersects the divisor at infinity exactly on the indetermination set of $f$, i.e. on the points where  $f$ cannot be holomorphically extended at infinity. This indetermination set is an analytical subset of $\mathbb{P}^n$ of (complex) dimension at most $n-2$.
(See Theorem 4.17 in \cite{G} for more details.)

The invariance property $G\circ f =\lambda_1 G$ implies that $f$ sends the sublevel set $\{G<\epsilon\}$ onto the sublevel set $\{G<\lambda_1 \epsilon\}$. In other words, the set $K^+:=\{G =0\}$, which we call {\it{filled Julia set}}, is included in an open set dilatated by  $f$. Therefore, every $y\notin K^+$ verifies $\displaystyle\lim_{m\rightarrow +\infty}\|f^m(y)\|=+\infty$. We thus define :
 
 \begin{defn}\label{Henon} A polynomial automorphism $f:\mathbb{C}^n\rightarrow \mathbb{C}^n$, $n\geq 2$, algebraically stable and such that $\lambda_1>\lambda_2$, will be called {\em{quasi-H\'enon}} if $X_f\cap I_f=\emptyset$ and the set $I_f$ is $f^{-1}$-attracting.
\end{defn}
 
\begin{exemple} The H\'enon's applications $f(x,y) = (x^2+ay, x)$, $a\in\mathbb{C}^*$ are quasi-H\'enon automorphisms of $\mathbb{C}^2$ if $|a| < 1$.\end{exemple} 

In the following subsections, $f:\mathbb{C}^n\rightarrow \mathbb{C}^n$ will be supposed to be a quasi-H\'enon map. We then have $\lambda_n=1$ and $f_*=(f^{-1})^*$ on $\mathbb{C}^n$.

\subsection{Equidistribution result for currents}\label{3.1}

\noindent The following result, which generalizes a theorem of \cite{BS} (see also Theorem 7.17 in \cite{FS}), will allow us to study locally the current $dd^c G$ in dimension $n> 1$.

\begin{lem}\label{eq2} Let $f:\mathbb{C}^n\rightarrow \mathbb{C}^n$ be a quasi-H\'enon map. Let $S$ be a closed positive $(1,1)$-current on $\mathbb{P}^n$, with mass 1, without Lelong numbers, and let $\xi\geq 0$ be a $C^{\infty}$ function of compact support in $\mathbb{C}^n\subset \mathbb{P}^n$. Then, denoting $\displaystyle c := \int_{\mathbb{P}^n}\xi S\wedge T^-\geq 0$, we have $$\displaystyle\lim_{m\rightarrow +\infty}\frac{1}{{\lambda_1}^m}f^{m*}(\xi S) =c \cdot T^+.$$
On the other hand, denoting $c' :=\displaystyle \int_{\mathbb{P}^n}\xi S\wedge T^+\geq 0$, we have $$\displaystyle\lim_{m\rightarrow +\infty}\frac{1}{{\lambda_{1}}^m}f_*^{m}(\xi T^-)=\lim_{m\rightarrow +\infty}\frac{1}{{\lambda_{1}}^m}f_*^{m}(\xi \omega^{n-l}) =c' \cdot T^-.$$   \end{lem}

\noindent \textit{Proof} We will give the proof in the case $n>2$, the case $n = 2$ being only a simplified version. We may suppose, without loss of generality, that $0< \xi\leq 1$. Let $\omega$ be the Fubiny-Study $(1,1)$-form. Let $\phi$ be a closed $(0,1)$-form.

We put $\displaystyle S_m:=\frac{1}{{\lambda_{1}}^m}f^{m*}(\xi S)$. Let us first show that the sequence of currents $(S_m)_m$ is bounded. It follows, using the fact that the currents $\displaystyle \frac{1}{{\lambda_1}^m}f^m_*\omega^{n-1}$ are cohomologous, that $\forall m\in\mathbb{N}$, $$\displaystyle \int_{\mathbb{P}^n} S_m\wedge \omega = \int_{\mathbb{P}^n}\xi S\wedge\frac{1}{{\lambda_1}^m}f^m_*\omega^{n-1}\leq \int_{\mathbb{P}^n} S\wedge\frac{1}{{\lambda_1}^m}f^m_*\omega^{n-1}   = \int_{\mathbb{P}^n} S\wedge \omega^{n-1} <+\infty.$$
The mass of the currents $S_m$ is thus uniformly bounded. 

Let us now show that an accumulation point $S_{\infty}$ of the sequence $(S_m)_m$ is closed. The Cauchy-Schwarz inequality gives : 
\begin{align*}
& \displaystyle\lim_{m\rightarrow +\infty}\bigg|\int_{\mathbb{P}^n}\partial S_m \wedge \phi \wedge \omega^{n-2}\bigg| = \lim_{m\rightarrow +\infty}\frac{1}{{\lambda_1}^m}\bigg|\int_{\mathbb{C}^n} \partial \xi\wedge S\wedge f^{m}_*(\omega^{n-2})\wedge f^{m}_* \phi\bigg|\\
& \displaystyle\leq \lim_{m\rightarrow +\infty}\frac{1}{{\lambda_1}^m}\bigg|\int_{\mathbb{C}^n}f^{m}_*(\omega^{n-2})\wedge S\wedge\partial\xi\wedge i\overline{\partial}\xi\bigg|^{1/2}\cdot\bigg|\int_{\mathbb{C}^n}f^{m}_*(\omega^{n-2}\wedge\phi\wedge i\overline{\phi})\wedge S\bigg|^{1/2}\\
& \displaystyle =\lim_{m\rightarrow +\infty}\bigg(\frac{\lambda'_{n-2}}{\lambda_1}\bigg)^{m/2}\bigg|\int_{\mathbb{C}^n}\bigg(\frac{1}{\lambda'_{n-2}}\bigg)^mf^{m}_*(\omega^{n-2})\wedge S\wedge\partial\xi\wedge i\overline{\partial}\xi\bigg|^{1/2}\\
&\ \ \ \ \ \ \ \ \ \ \ \ \cdot\bigg|\int_{\mathbb{C}^n}\bigg(\frac{1}{\lambda_1}\bigg)^mf^{m}_{*}\big(\omega^{n-2}\wedge\phi\wedge i\overline{\phi}\big)\wedge S\bigg|^{1/2}\\
& \displaystyle =\lim_{m\rightarrow +\infty}\bigg(\frac{\lambda_{2}}{\lambda_1}\bigg)^{m/2}\bigg|\int_{\mathbb{C}^n}\bigg(\frac{1}{\lambda_{2}}\bigg)^mf^{m}_*(\omega^{n-2})\wedge S\wedge\partial\xi\wedge i\overline{\partial}\xi\bigg|^{1/2}\\
& \ \ \ \ \ \ \ \ \ \ \ \ \cdot\bigg|\int_{\mathbb{C}^n} \omega^{n-2}\wedge\phi\wedge i\overline{\phi}\wedge\bigg(\frac{1}{\lambda_1}\bigg)^mf^{m*} S\bigg|^{1/2}\\
& \displaystyle =\lim_{m\rightarrow +\infty} O\bigg(\bigg(\frac{\lambda_{2}}{\lambda_1}\bigg)^{m/2}\bigg) =0.
\end{align*}

We will further show that $S_{\infty}=c\cdot T^+$. By virtue of Theorem 1.1 in \cite{B}, we have $\displaystyle \frac{1}{{\lambda_1}^m}f^{m*}S\longrightarrow T^+$ in the sense of currents. Since $\displaystyle S_m\leq \frac{1}{{\lambda_1}^m}f^{m*}S$ for all $m$, due to the positivity of the operator $f^*$, we have $S_{\infty}\leq T^+$. The extremality of $T^+$ (Theorem 2.1 in \cite{G2}) leads to $S_{\infty}=c\cdot T^+$, with $c: = \int_{\mathbb{P}^n}\xi S\wedge T^-\geq 0$. (Note that for $n=2$, the assumption that $S$ does not have Lelong numbers can be removed thanks to Theorem A in \cite{P2}).

In order to show that $\displaystyle\lim_{m\rightarrow +\infty}\frac{1}{{\lambda_1}^m}{f_*}^{m}(\xi S) =c \cdot T^-$ when $S=\omega$ or $T^-$, the proof proceeds similarly, with the following precision. If $S_{\infty}$ is an accumulation point of the sequence of term $\frac{1}{{\lambda_1}^m}{f_*}^{m}(\xi S)$, then $S_{\infty}=c'\cdot T^-$ follows from Theorem 2.6 in \cite{G2}.$\blacksquare$\\ 

\noindent Let us give some consequences. Recall that we denote $K^+:=\{G=0\}$.

\begin{pro}\label{pro}
Let $f:\mathbb{C}^n\rightarrow \mathbb{C}^n$, $n\geq 1$ be a quasi-H\'enon map. Then $\partial K^+$ is the support of $T^+$. 
\end{pro}

\noindent \textit{Proof} Recall that we denote $J:=supp(T^+)=supp(dd^c G)$. Note that the sets $K^+$, $\partial K^+$ and $J$ are invariant by $f$ and by $f^{-1}$. 

There exists a $C^{\infty}$ function $\xi\geq 0$ with compact support in $\mathring{K}^+$ whose interior intersects $J^-$. Using Lemma \ref{eq2}, if $\omega$ is the Fubiny-Study $(1,1)$-form, there exists $c>0$ such that : $$\lim_{m\rightarrow +\infty}\frac{1}{{\lambda_{1}}^m}f^{m*}(\xi \omega) = c\cdot T^+.$$ It follows that $J\subset K^+$, and more precisely that $J\subset\partial K^+$.\\

\noindent Conversely, suppose by contradiction that there exists $x\in \partial K^+\setminus J$. The function $G$ is harmonic in a small neighborhood $W$ of $x$, where it has a minimum $G(x)=0$, hence it vanishes in $W$, which is impossible. Thus $J=\partial K^+$.$\blacksquare$\\ 

\begin{lem}\label{julia}
The filled Julia set $K^+$ of a quasi-H\'enon map $f:\mathbb{C}^n \to \mathbb{C}^n  $, $n\geq 2$, if it is of non empty interior, satisfies $\overline{\mathring{K}^+}=K^+$.
\end{lem}

\noindent \textit{Proof} Let $\xi\omega$ be a $(1,1)$-form of support in $\mathring{K}^+$, where $\omega$ is the Fubini-Study form and  $\xi\geq 0$ is a $C^{\infty}$ function of compact support in $\mathring{K}^+$ intersecting $J^-$. Note that $K^+$ is invariant by $f$ and by $f^{-1}$, thus $\mathring{K}^+$ is also. Note that, by Proposition \ref{pro}, $J=\partial K^+$. The support of $f^{m*}(\xi\omega)$ being in $\mathring{K}^+$ for all $m$, by Lemma \ref{eq2}, every open subset of  $\mathbb{C}^n$  intersecting $J=\partial K^+$ also  intersects $\mathring{K}^+$, thus $\overline{\mathring{K}^+}=K^+$.$\blacksquare$\\

\subsection{Equidistribution result for measures}\label{3.2}

We give conditions under which the means of pull-backs $\displaystyle \mu_m:=\frac{1}{m}\sum_{i=0}^m f^{i*} \mu$ converge to the invariant measure  $T^+\wedge T^-$. Recall that we denote by $J$ and respectively $J^-$ the supports of the currents $T^+$ and $T^-$.

\begin{lem}\label{dense} Let $f:\mathbb{C}^n\rightarrow \mathbb{C}^n$ be a quasi-H\'enon map. Let $U$ be an open set  intersecting $J\cap J^-$. Then $\displaystyle\bigcup_{m\geq 0}f^m(U|_{J^-})$, and even $\displaystyle\bigcap_{i\geq 0}\bigcup_{m\geq i}f^m(U|_{J^-})$, are dense in $J^-$ for the induced topology.\\
More precisely, for every open set $V$ intersecting $J^-$, there exists $\epsilon>0$ such that $$\displaystyle \|T^-\|\left(\displaystyle f^{m}(U|_{J^-})\cap V\right)>\epsilon$$ for $m$ sufficiently large.\end{lem}

\noindent \textit{Proof} Let $\xi$ be a positive  $\mathcal{C}^{\infty}$ function with compact support $Z$, such that $U:=\mathring{Z}$ intersects $supp(T^+\wedge T^-)\subset J\cap J^-$. By Lemma \ref{eq2}, the sequence of currents $\displaystyle\left(\frac{1}{{\lambda_{1}}^m}{f^{m}}_*\left(\xi T^-\right)\right)_m$ converges to $c\cdot T^-$, with $c = \int \xi T^-\wedge T^+>0$, when $m\rightarrow +\infty$. 

But the support of $\displaystyle\frac{1}{{\lambda_{1}}^m}{f^{m}}_*(\xi T^-)$ is $f^{m}(Z)\cap J^-$, hence for every open set $V$ intersecting $J^-$, $\exists m_0\in \mathbb{N}$ such that $f^{m}(U|_{J^-})\cap V\neq \emptyset$ when $m\geq m_0$. In other words, $f^{m}(U|_{J^-})\cap V|_{J^-}\neq \emptyset$ when $m\geq m_0$. Therefore, $\displaystyle\bigcup_{m\geq 0}f^m(U|_{J^-})$ is dense in $J^-$ for the induced topology, and even $\displaystyle\bigcap_{i\geq 0}\bigcup_{m\geq i}f^m(U|_{J^-})$.

More precisely, it follows from the  convergence of $\displaystyle\left(\frac{1}{{\lambda_{1}}^m}{f^{m}}_*\left(\xi T^-\right)\right)_m$ to $c\cdot T^-$ when $m\rightarrow +\infty$ that there exists $\epsilon>0$ such that $\displaystyle\|T^-\|\left(\displaystyle f^{m}(U|_{J^-})\cap V\right)>\epsilon$ when $m\geq m_0$.\\
Indeed, let us suppose that $\displaystyle \max_{x\in \mathbb{C}^n}\xi(x) =1$, and consider a test function $\rho$ positive on $V$, such that $\displaystyle \max_{x\in \mathbb{C}^n}\rho(x)=1$. Then : 
\begin{align*}
\displaystyle \|T^-\|\left(\displaystyle f^{m}(U|_{J^-})\cap V\right)&= \frac{1}{{\lambda_{1}}^m}\int_{f^m(U)\cap V}({{f^{m}}_*} T^-)\wedge \omega\geq \frac{1}{{\lambda_{1}}^m}\int_{\mathbb{C}^n}({{f^{m}}_*} \xi T^-)\wedge \rho\omega\\
& \displaystyle\underset{m\rightarrow +\infty}{\longrightarrow} c\int_{\mathbb{C}^n}T^-\wedge\rho\omega>0 .
\end{align*}
$\blacksquare$\\ 

\noindent We can now give the following convergence result (see for example \cite{W} or \cite{M} for the definition of an ergodic measure). The fact that the support of $T^+\wedge T^-$ is $J\cap J^-$ is proved in dimension $2$ in \cite{F} in a little different context.

\begin{pro}\label{supportmesure} Let $f:\mathbb{C}^n\rightarrow \mathbb{C}^n$ be a quasi-H\'enon map. If $\mu$ is a probability measure with compact support on $J^-$ that does not charge pluripolar sets, then the sequence of measures $\displaystyle \left(\mu_m:=\frac{1}{m}\sum_{i=0}^m f^{i*} \mu\right)_m$ converges weakly to the invariant probability measure $T^+\wedge T^-$.\\
Moreover, the measure $T^+\wedge T^-$ is ergodic, and its support is the compact set $J\cap J^-$.\end{pro}

\noindent \textit{Proof} Let us denote by $\mu_f$ an accumulation point of the sequence of probability measures $\displaystyle \left(\frac{1}{m}\sum_{i=0}^m f^{i*} \mu\right)_m$. By construction, we have $f^*\mu_f=\mu_f$, with $supp(\mu_f)\subset J\cap J^-$. It remains to show the inclusion in the opposite direction.\\

\noindent 1) {\it{Let us show that the support of $T^+\wedge T^-$ is exactly $J\cap J^-$.}} Let $U\subset \mathbb{C}^n$ be an open set intersecting $J\cap J^-$. Denote by $\xi$ a positive test function such that $\xi\equiv 1$ on $U$. By Lemma \ref{eq2}, we have  $$\displaystyle  \frac{1}{{\lambda_{1}}^m}{f^m}_*(\xi T^-)\rightarrow c \cdot  T^-$$ for a constant $c\geq 0$. We will show that $c>0$. Since $\displaystyle c=\int_{\mathbb{C}^n}\xi d(T^+\wedge T^-)$, this will imply that the support of $T^+\wedge T^-$ is included in $J\cap J^-$.\\
Note that  $\forall m\in \mathbb{N}$, we have $$\displaystyle \frac{1}{{\lambda_{1}}^m}{f^m}_*(\xi T^-)\geq \frac{1}{{\lambda_{1}}^m}{f^m}_*T^-|_U=T^-|_{f^m(U)}.$$ But by Lemma \ref{dense}, for every open set $V$ intersecting $J^-$, $\exists \epsilon >0$, such that for $m$ sufficiently large, we have \begin{equation}\label{epsilon}\displaystyle \|T^-\|\left(\displaystyle f^{m}(U|_{J^-})\cap V\right)>\epsilon.\end{equation}Thus $c>0$. Finally $supp(T^+\wedge T^-)\subset J\cap J^-$. The inverse inclusion is trivial.\\

\noindent 2) Let $\rho$ be a test function whose support $U$ intersects $J^+\cap J^-$. We can decompose $\rho$ into the difference of two positive continuous plurisubharmonic functions $\rho=u-v$ (for example, we can choose $u=\rho+A\log(1+\|z\|^2)$ and $v = A\log(1+\|z\|^2)$ for $A$ sufficiently large). Consider the sequence of functions $\displaystyle(\phi_m)_m:=\left(\frac{1}{m}\sum_{i=1}^m {f^i}_* u|_{ J^-}\right)_m$ defined on $\mathcal{S}$. \\

Note that this sequence is bounded on every compact of $J^-$ : indeed, let $U$ be an open set containing $J\cap J^-$, and denote $H$ the compact set $J^-\cap \{G=0\}$. Then $U|_{J^-}\subset\subset f(U)|_{J^-}\cup H$. Consequently, we have $\forall m\in\mathbb{N}$, $\displaystyle \sup_{U|_{J^-}}\phi_m \leq \sup_{U|_{J^-}\cup H}u$.\\
Denote $\displaystyle{\phi_{\infty}}:={\limsup_{m\rightarrow +\infty}}\text{ }\phi_m$. Then, by Lemma \ref{sup}, $\phi_{\infty}^*$ is a $T$-psh function. Since $\phi_{\infty}$ satisfies by construction $\displaystyle f_*\phi_{\infty}=\phi_{\infty}$, we have by Lemma \ref{sup} : \begin{equation}\label{invariant} f_*\phi_{\infty}^*=\phi_{\infty}^*.\end{equation} 

\noindent 3) {\it{Let us show that $\phi_{\infty}^*$ is constant on $J^-$.}} Let $\displaystyle c:=\max_{J\cap J^-}\phi_{\infty}^*$ and $\epsilon>0$. Then $J\cap J^-$ is included in the open set $U:=\{\phi^*_{\infty}<c+\epsilon\}$. Then (\ref{invariant}) implies that $\forall m\in\mathbb{N}$, $$f^m(U)=\{{f^m}_*\phi^*_{\infty}<c+\epsilon\}=U.$$
We have $\phi^*_{\infty}\leq c+\epsilon$ on $\displaystyle U=\bigcup_{m\geq 0}f^m(U)$, and this open set being dense in $J^-$ for the induced  topology by  Lemma \ref{dense}, Theorem \ref{hartogs} and Corollary \ref{hartogs2} imply $\displaystyle\phi^*_{\infty}\leq c+\epsilon$ on $J^-$. The constant $\epsilon>0$ being arbitrary, it follows that $\phi^*_{\infty}\leq c$ on $J^-$. Instead of invoking Proposition \ref{liouville}, let us pursue with some dynamical arguments.\\

\noindent Let $d:=\displaystyle\inf_{J\cap J^-}\phi^*_{\infty}$ and $\epsilon>0$. Let $V:=\{\phi^*_{\infty}< d-\epsilon\}$. Then $J\cap J^-\cap\overline{ V}=\emptyset$. As previously we have $V= f(V)$; since for any bounded set $W$ containing $J\cap J^-$ there exists $m_0\in\mathbb{N}$ such that $f^m(V)\cap W\neq\emptyset$ for $m\geq m_0$, we obtain $V=\emptyset$. Therefore, $\phi^*_{\infty}\geq d$ on $J^-$.\\

\noindent If $\phi^*_{\infty}|_{J\cap J^-}$ is not constant on $J\cap J^-$, there exists $\epsilon>0$ such that the open set $W:=\{\phi^*_{\infty}< c-\epsilon\}$ is non empty and intersects $J\cap J^-$. Moreover, we have again $f(W)=W$. By Lemma \ref{dense}, $\displaystyle W=\bigcup_{m\geq 0}f^m(W)$ is an open set dense in $J^-$. By  Theorem \ref{hartogs} and Corollary \ref{hartogs2}, $\phi^*_{\infty}<c-\epsilon$ in $J\cap J^-$, which contradicts the definition of $\displaystyle c:=\max_{J\cap J^-}\phi^*_{\infty}$. Therefore, $\phi^*_{\infty}=c$ on $J\cap J^-$. Thus $c=d$, hence $\phi^*_{\infty}=c$ on $J^-$.\\

\noindent 4) {\it{Conclusion.}} The sequence $(\phi_m)_m$ converges in $L^1(T^+\wedge T^-)$ by Birkhoff's ergodic theorem. The set $\{\phi_{\infty}<\phi^*_{\infty}\}$ is pluripolar by Lemma \ref{trucdeouf}. Since the measure $T^+\wedge T^-$ does not charge pluripolar sets, we deduce that the sequence $(\phi_m)_m$ converges $T^+\wedge T^-$-a.e. towards $\displaystyle c={\limsup_{m\rightarrow +\infty}}^*\phi_m$. Proposition \ref{liouville} implies that this sequence converges $\|T^-\|$-a.e. towards $c$ on $J^-$. Let us denote by $(\mu_{m_i})_i$ a sequence of measures, extracted from $(\mu_m)_m$, which converges to $\mu_f$. Since the measure $\mu$ does not  charge, by hypothesis, the pluripolar sets, thanks to the dominated convergence theorem, there exists a constant $C\in\mathbb{R}$ such that $$\displaystyle <\mu_f,\rho> = \lim_{i\rightarrow +\infty} < \mu,\frac{1}{m_i}\sum_{j=1}^{m_i} {{f^{j}}_*} u>- \lim_{i\rightarrow +\infty} < \mu,\frac{1}{m_i}\sum_{j=1}^{m_i} {{f^{j}}_*} v>= C.$$
On the other hand, since the measure $T^+\wedge T^-$ does not charge pluripolar sets (see e.g. \cite{G}), we have $$\displaystyle <T^+\wedge T^-,\rho> = \lim_{i\rightarrow +\infty} < T^+\wedge T^-,\frac{1}{m_i}\sum_{j=1}^{m_i} {{f^{j}}_*} u>- \lim_{i\rightarrow +\infty} < T^+\wedge T^-,\frac{1}{m_i}\sum_{j=1}^{m_i} {{f^{j}}_*} v>= C.$$ Thus $T^+\wedge T^-=\mu_f$. The ergodicity of this invariant measure follows directly from Theorem 4.4.1 in \cite{M}.$\blacksquare$\\ 

\begin{rmq}\label{towards} Let us emphasize that, in the proof of Proposition \ref{supportmesure}, we have obtained that if $\rho$ is a test function, the sequence of functions $\displaystyle\left(\frac{1}{m}\sum_{i=1}^m {{f^{i}}_*\rho}\right)_m$ converges towards a constant $C$ outside a pluripolar set, and therefore converges in $L^{1}_{loc}(T^+\wedge T^-$) towards $C$ thanks to the dominated convergence theorem. Moreover, $C>0$ if and only if the support of $\rho>0$ intersects $J\cap J^-$.\\
We thus deduce an alternative proof of the point 1) of the proof of Proposition \ref{supportmesure} above. Let $\xi>0$ be a test function with support intersecting $J\cap J^-$. Then $\forall m\in\mathbb{N}$, $$\displaystyle \int_{\mathbb{C}^n} \xi d(T^+\wedge T^-)=\int_{\mathbb{C}^n}\frac{1}{m}\sum_{i=1}^m  {f^i}_*\xi d(T^+\wedge T^-).$$
By taking $m\rightarrow +\infty$, since $T^+\wedge T^-$ does not charge pluripolar sets, the dominated convergence theorem implies that $\displaystyle \int_{\mathbb{C}^n} \xi d(T^+\wedge T^-)>0$, which concludes.\end{rmq}

\noindent The following remark generalizes Theorem 6.6 in \cite{FS2}. 
\begin{rmq}\label{deduction} We deduce from Remark \ref{towards} that for every open set $U$ intersecting $J\cap J^-$, the set $\displaystyle J^-\setminus \bigcup_{m\geq 0}f^m(U|_{J^-})$ is pluripolar. 
\end{rmq}

\noindent \textit{Proof} Indeed, let $\rho$ be a strictly positive test function of support included in $U$, and whose interior intersects $J\cap J^-$. Then
$$\displaystyle \phi_{\infty}:=\limsup_{m\rightarrow +\infty}\frac{1}{m}\sum_{i=1}^m  {{f^{i}}_*\rho}\leq \limsup_{m\rightarrow +\infty}\frac{1}{m}\sum_{i=1}^m  {{f^{i}}_*\mathbf{1}_U}.$$
On the other hand, Remark \ref{towards} shows that $\phi_{\infty}$ is equal, outside a pluripolar set, to a constant $C> 0$. Thus $$\displaystyle  J^-\setminus \bigcup_{m\geq 0}f^m(U|_{J^-})\subset
 \left\lbrace\limsup_{m\rightarrow +\infty}\frac{1}{m} \sum_{i=1}^m {{f^{i}}_*\mathbf{1}_U}=0\right\rbrace\subset  
\left\lbrace\phi_{\infty}=0\right\rbrace\subset 
\left\lbrace\phi_{\infty}<C\right\rbrace.$$
Thus all these sets are pluripolar, which proves the statement.$\blacksquare$\\




\begin{thebibliography}{}

\bibitem{A}{H. Alexander}, Projective capacity, Ann. Math. Stud., 100, 3--27 (1981).
\bibitem{B}{T. Bayraktar}, Equidistribution towards the Green current in big cohomology classes, Int. J. Math., 24, 1350080 (2013). 
\bibitem{Be}{A.F. Beardon}, Iteration of rational functions. Springer (1991).
\bibitem{BG2}{M.T. Belghiti}, \'El\'ements pour une th\'eorie constructive des fonctions lisses, Th\`ese de doctorat d'\'Etat, Universit\'e Ibn Tofa\"il, Kenitra (2004).
\bibitem{BG}{L. Gendre}, In\'egalit\'es de Markov singuli\`eres et approximation des fonctions holomorphes de la classe M, PhD thesis, Universit\'e Toulouse III (2005).
\bibitem{BK}{L. Bia{\l}as-Cie\.z, M. Kosek}, Iterated function systems and {\L}ojasiewicz-Siciak condition of Green's function, Potential Analysis, 34, 207--221 (2011).
\bibitem{BS}{E. Bedford, J. Smillie}, Polynomial diffeomorphisms of $\mathbb{C}^2$. II, J. Amer. Math. Soc., 4, 657--679 (1991).
\bibitem{Bo}{B. Berndtsson, N. Sibony}, The $\overline{\partial}$-equation on a positive current, Invent. Math., 147, 371--428 (2002).
\bibitem{D}{J.P. Demailly}, Monge-Ampere operators, Lelong numbers, and intersection theory, {\it{in}} Complex analysis and geometry, 115--191. V. Ancona, A. Silva ed., Plenum Press (1993).
\bibitem{DDS}{T.C. Dinh, Dujardin R., N. Sibony}, On the dynamics near infinity of some polynomial mappings in $\mathbb{C}^2$, Math. Ann., 333, 703--739 (2005).
\bibitem{DL}{T.C. Dinh, M.G. Lawrence}, Polynomial hulls and positive currents, Ann. Fac. Sci. Toulouse Math. series 6, 12, 317--334 (2003).
\bibitem{DS76}{T.C. Dinh, N. Sibony}, Green currents for holomorphic automorphisms of compact K\"ahler manifolds, J. Amer. Math. Soc., 18, 291--312 (2005).

\bibitem{DS4}{T.C. Dinh, N. Sibony}, Dynamique des applications polynomiales semi-r\'eguli\`eres, Ark. Mat., 42, 61--85 (2004).
\bibitem{DS2}{J. Duval, N. Sibony}, Polynomial convexity, rational convexity, and currents, Duke Math. J., 79, 487--513 (1995).
\bibitem{E}{D.A. Edwards}, Choquet boundary theory for certain spaces of lower semicontinuous functions, \textit{in} Function algebras, Birtel ed., Scott Foresman and Cie, 300--309 (1966).
\bibitem{F}{J.E. Fornaess}, The Julia set of H\'enon maps, Math. Ann., 334, 457--464 (2006).
\bibitem{FaG}{C. Favre, V. Guedj}, Dynamique des applications rationnelles des espaces multiprojectifs, Indiana Univ. Math. J., 50, 881--934 (2001).
\bibitem{FS}{E. Fornaess, N. Sibony}, Complex dynamics in higher dimensions, \textit{in}  Complex potential theory, NATO ASI series C, 439, 131--18. Gauthier ed. 6 (1994).
\bibitem{FS2}{E. Fornaess, N. Sibony}, Complex dynamics in Higher dimensions II, \textit{in} Modern methods in complex analysis, 135--182. Annals of mathematics studies, 137 (1995).
\bibitem{G}{V. Guedj}, Propri\'et\'es ergodiques des applications rationnelles, \textit{in} Quelques aspects des syst\`emes dynamiques polynomiaux, Panor. Synth\`eses, 30, 97--202. Soc. Math. France, Paris (2010).
\bibitem{GS2}{V. Guedj, N. Sibony}, Dynamics of polynomials automorphisms of $\mathbb{C}^k$, Ark. Mat., 40, 207--243 (2002).
\bibitem{G2}{V. Guedj}, Courants extr\'emaux et dynamique complexe, Ann. Sci. \'Ecole Norm. Sup., 38, 407--426 (2005).
\bibitem{GR}{T.W. Gamelin, H. Rossi}, Jensen measures and algebras of analytic functions, \textit{in} Function algebras, 15--35. F.T. Birtel ed., Scott Foresman and Cie (1966).
\bibitem{K}{M. Klimek}, Pluripotential theory. Clarendon Press (1991).
\bibitem{M}{A. Lasota, M.C. Mackey}, Probabilistic properties of deterministic systems. Cambridge University Press (1985).
\bibitem{Po}{E.A. Poletsky}, Plurisubharmonic functions as solutions of variational problems, \textit{in} Several complex variables and complex geometry (Santa Cruz, CA, 1989), Part. 1,  163--171. American. Math. Soc. (1991).
\bibitem{Po2}{E.A. Poletsky}, Holomorphic currents, Indiana Univ. Math. J., 42, 85--144 (1993).
\bibitem{P}{F. Protin}, Dynamique d'endomorphismes polynomiaux de $\mathbb{C}^k$, PhD thesis, Universit\'e Toulouse III (2010).
\bibitem{P2}{F. Protin}, Equidistribution vers le courant de Green, Ann. Polon. Math., 115, 201--218 (2015).
\bibitem{Ra}{T.J. Ransford}, Jensen measures, \textit{in} Approximation, complex analysis and potential theory, 221--237. N. Arakelian and P.M. Gauthier ed., Kluwer Academic Publishers (2001).
\bibitem{Ro}{J.P. Rosay}, Poletsky theory of disks on holomorphic manifolds, Indiana Univ. Math. J., 52, 157--169 (2003).
\bibitem{S2}{N. Sibony}, arXiv:1509.01790
\bibitem{St}{E.L. Stout}, Polynomial convexity. Birka\"user Basel, Progress in mathematics, 261 (2007).
\bibitem{St2}{E.L. Stout}, The theory of uniform algebras. Bogden and Quigley (1971).
\bibitem{W}{Walters}, An introduction to ergodic theory. Springer (1981).

\end{thebibliography}
\end{document}